\documentclass[twoside,11pt,a4wide]{article}
\usepackage{amsfonts}
\usepackage{amssymb}
\usepackage{a4wide}
\usepackage{amsthm}
\usepackage{amsmath}
\usepackage[all]{xy}
\usepackage{enumerate}
\usepackage{fancyhdr}

\newtheorem{theorem}{Theorem}[section]
\newtheorem{proposition}[theorem]{Proposition}
\newtheorem{corollary}[theorem]{Corollary}
\newtheorem{lemma}[theorem]{Lemma}
\theoremstyle{definition}
\newtheorem*{notation}{Notation}
\newtheorem*{Beweis}{Proof}
\newtheorem{definition}[theorem]{Definition}
\newtheorem{punto}[theorem]{}
\theoremstyle{remark}
\newtheorem{remark}[theorem]{Remark}
\newtheorem{ex}[theorem]{Example}
\newtheorem{exs}[theorem]{Examples}
\newtheorem{remarks}[theorem]{Remarks}
\CompileMatrices
\swapnumbers
\CompileMatrices
\input{tcilatex}
\begin{document}

\title{Hopf Semialgebras\thanks{%
MSC2010: primary 16T05 (secondary 16T10, 16T15, 16T20)}}
\author{\textbf{Jawad Y. Abuhlail}\thanks{%
Corresponding Author} \ \thanks{%
The author would like to acknowledge the support provided by the Deanship of
Scientific Research (DSR) at King Fahd University of Petroleum $\&$ Minerals
(KFUPM) for funding this work through project No. IN080400.} \\
%EndAName
Department of Mathematics and Statistics\\
Box $\#\;$5046, KFUPM; 31261 Dhahran, KSA \\
abuhlail@kfupm.edu.sa \\
\\
\textbf{Nabeela Al-Sulaiman} \\
Department of Mathematics\\
University of Dammam; 31451 Dammam, KSA \\
nalsulaiman@ud.edu.sa}
\date{\today }
\maketitle

\begin{abstract}
In this paper, we introduce and investigate \emph{bisemialgebras }and\emph{\
Hopf semialgebras} over commutative semirings. We generalize to the
semialgebraic context several results on bialgebras and Hopf algebras over
rings including the main reconstruction theorems and the \emph{Fundamental
Theorem of Hopf Algebras}. We also provide a notion of \emph{quantum monoids}
as Hopf semialgebras which are neither commutative nor cocommutative; this
extends the Hopf algebraic notion of a quantum group. The generalization to
the semialgebraic context is neither trivial nor straightforward due to the
non-additive nature of the base category of Abelian monoids which is also
neither Puppe-exact nor homological and does not necessarily have enough
injectives.
\end{abstract}

\section*{Introduction}

\qquad Topological investigations of Lie groups and group-like spaces led
the German mathematician Heinz Hopf to realize that the multiplication map
in the cohomology algebra yields a comultiplication, and from that
combination Hopf got remarkable structure results \cite{Hop1941}. This was
the formal birth of the theory of \emph{Hopf algebras}, one of the main
streams of research in mathematics nowadays. Apart from their nice theory
from the purely algebraic point of view \cite{Swe1969}, \cite{Rad2012} (%
\emph{e.g. }Kaplansky's conjectures, Andruskiewitsch-Schneider's project on
the classification of semisimple finite dimensional Hopf algebras, the
module theoretic approach \cite{BW2003}), Hopf algebras play important roles
in many aspect of mathematics like graded ring theory (coactions \cite%
{Mon1993}, \cite{DNR2001}), algebraic geometry (affine group schemes \cite%
{Abe1980}, \cite{Und2011}), number theory (formal groups), mathematical
physics (quantum groups \cite{Maj1990}), Lie algebras (universal enveloping
algebras), Topology (\emph{e.g.} cohomology of exceptional Lie groups),\
Knot Theory \cite{KRT1997}, non-commutative geometry, Galois theory
(Hopf-Galois extensions), combinatorics (umbral calculus), computer science (%
\emph{e.g.} models of linear logic) and many more \cite{CGW2006}.

Moreover, in category theory, \emph{Hopf monoids} in braided monoidal
categories \cite{Tak2000} and \emph{Hopf monads} in arbitrary categories 
\cite{MW2011} are gaining increasing interest \cite{Ver}. While some basic
definitions and results remain the same when moving from one category to
another \cite{Tak1999}, \cite{CGW2006}, several structural properties of the
category of Hopf monoids depend naturally on the properties of the category
in which such objects live (\emph{e.g.} the category $\mathbf{Hopf}_{R}$ of
Hopf algebras over a commutative ring $R$ is closed under limits in the
category $\mathbf{Bilag}_{R}$ of $R$-bialgebras if $R$ is von Neumann
regular \cite{Por2011a}).

In this paper, we introduce and investigate \emph{Hopf semialgebras} (\emph{%
bisemialgebras}) over commutative semirings. Let $S$ be a commutative
semiring and denote by $\mathbb{S}_{S}$ the category of $S$-semimodules and
by $\mathbb{CS}_{S}\overset{\iota }{\hookrightarrow }\mathbb{S}_{S}$ the 
\emph{full} subcategory of cancellative $S$-semimodules \cite{Tak1981}. A
notion of Hopf semialgebras (bisemialgebras) was introduced by the first
author\footnote{%
http://www.ingvet.kau.se/juerfuch/conf/nomap/talk/Abuhlail.pdf} using
Takahashi's tensor-like product \cite{Tak1982}, which we denote in this
paper by $\boxtimes _{S}.$ That notion was investigated by the second author
in her dissertation \cite{AlS2011} in the category $\mathbb{CS}_{S}$
assuming also that the base semiring $S$ is \emph{cancellative.} In this
paper, we use instead the natural tensor product $\otimes _{S}$ in $\mathbb{S%
}_{S}$ inherited from the tensor product in the symmetric monoidal category $%
(\mathbf{AbMonoid},\otimes ,\mathbb{N}_{0})$ of Abelian monoids \cite%
{Kat1997}. As clarified in \cite{Abu-a}, $-\boxtimes _{S}-$ and $-\otimes
_{S}-$ are isomorphic bifunctors on $\mathbb{CS}_{S},$ whence the results in
this paper generalize those in \cite{AlS2011}. A main advantage of using $%
\otimes _{S}$ (instead of $\boxtimes _{S}$) is that the category $(\mathbb{S}%
_{S},\otimes _{S},S)$ is a symmetric monoidal category (while the category $%
(S_{S},\boxtimes _{S},S)$ is, in general, \emph{semiunital semimonoidal} 
\cite{Abu2013}). This suggests defining Hopf semialgebras (bisemialgebras)
as Hopf monoids (bimonoids) in $(\mathbb{S}_{S},\otimes _{S},S).$

Such Hopf monoids (bimonoids) not only add \emph{new families of concrete
examples} to the literature, but they are of particular importance for
theoretical and practical reasons. On one hand, and in contrast to the
category of modules over a ring, the category of semimodules over a semiring
is not Abelian (not even additive) and so many proofs that depend heavily on
lemmas of diagrams cannot be directly applied to our context. Add to that
this category is not \emph{Puppe-exact} and not \emph{homological} \cite%
{BB2004} and so a new notion of exact sequences for semimodules over a
semiring was necessary to prove restricted versions of the Short Five Lemma
and the Snake Lemma \cite{Abu-b}. Moreover, working over proper semifields
(semisimple proper semirings) does not bring big advantages as was the case
in the theory of Hopf algebras over fields (semisimple rings). This is due
to the fact that all semivector spaces over a semifield $F$ (semimodules
over a semisimple semiring $S$) are free (projective) if and only if $F$ is
a field \cite[Theorem 5.11]{KN2011} ($S$ is a semisimple ring \cite[Theorem
5.7]{KN2011}). This suggests that one uses a combination of techniques from
categorical algebra, universal algebra and homological algebra to overcome
these and other difficulties \cite{Abu-c}. On the other hand, semirings and
semimodules proved to have a wide spectrum of significant applications in
several aspects of mathematics like optimization theory, tropical geometry,
idempotent analysis, physics, theoretical computer science (\emph{e.g.}
Automata Theory) and many more \cite{Gol1999}. It is hoped that
investigating Hopf semialgebras would bring new applications (for some
applications of bisemialgebras in Automata Theory, see \cite{Wor2012}).

This paper is organized as follows: after this introduction, we recall in
Section 1 some basic definitions and properties of semicoalgebras and
semicomodules (for a detailed discussion of semicorings and semicomodules,
see \cite{Abu-c}). In Section 2, we introduce the notion of a bisemialgebra
and study \emph{integrals} \textbf{on} and \textbf{in} a given
bisemialgebra. Moreover, we give a reconstruction result for bisemialgebras
using the notion of a \emph{bimonad }in the sense of \cite{BBW2009} (Theorem %
\ref{bimonad}). In Section 3, we investigate the categories of Doi-Koppinen
semimodules and give relatively weak sufficient conditions for such a
category to be a Wisbauer category of type $\sigma \lbrack M]$ for a
suitable \emph{subgenerator} $M$ (Theorem \ref{AC-sigma}). In Section 4, we
consider Hopf semialgebras and extend several examples of quantum groups to 
\emph{quantum monoids} which we introduce as non-commutative
non-cocommutative Hopf semialgebras. Moreover, we present the \emph{%
Fundamental Theorem of Hopf Semialgebras }(Theorem \ref{FTHA}).\emph{\ }A
reconstruction result for Hopf semialgebras in terms of Hopf monads \cite%
{MW2011} is also obtained (Theorem \ref{H-monad}). In addition to that, we
use integrals to characterize Hopf semialgebras which are cosemisimple as
semicoalgebras (Proposition \ref{H-cosemi}) and those which are semisimple
as semialgebras (Proposition \ref{H-semi}). In Section 5, we present
possible constructions of \emph{dual} semicoalgebras, dual bisemialgebras
and dual Hopf semialgebras in both the finite and the infinite cases.

\section{Preliminaries}

\qquad In this section, we collect some definitions and properties of
semirings (semimodules) and semicoalgebras (semicomodules).

\subsection*{Semirings and Semimodules}

\qquad A \emph{semiring} is essentially a \emph{monoid} in the category $%
\mathbf{AbMonoid}$ of Abelian Monoids, or -- roughly speaking -- a ring not
necessarily with subtraction. Moreover, a \emph{semifield} is a semiring in
which every non-zero element is invertible. Indeed, any ring is a semiring;
the first natural examples of semirings (semifields) which are not rings
(not fields) are the set $\mathbb{N}_{0}$ of natural numbers ($\mathbb{R}%
^{+}:=[0,\infty )$ and $\mathbb{Q}^{+}:=\mathbb{Q}\cap \lbrack 0,\infty )$)
with the usual addition and multiplication. All semirings in this paper are
unital, and for a given semiring $S,$ we assume that $0_{S}\neq 1_{S}.$
Given a semiring $S$, we mean by a right (left) $S$-semimodule a right
(left) $S$-module not necessarily with subtraction. The category of right
(left) $S$-semimodules is denoted by $\mathbb{S}_{S}$ ($_{S}\mathbb{S}$).
For two semirings $S$ and $T,$ an $(S,T)$-bisemimodule is a left $S$%
-semimodule, which is also a right $T$-semimodule such that $(sm)t=s(mt)$
for all $s\in S,$ $m\in M$ and $t\in T.$ The category of $(S,T)$%
-bisemimodules and $S$-linear $T$-linear maps (called $(S,T)$\emph{-bilinear
maps}) is denoted by $_{S}\mathbb{S}_{T}.$ We refer the reader to \cite%
{Gol1999} and the first section of \cite{Abu-c} for the basic definitions
and properties of semirings and semimodules. For any right (left) $S$%
-semimodule, we have a canonical isomorphism of Abelian monoids $M\otimes
_{S}S\overset{\vartheta _{M}^{r}}{\simeq }M$ ($S\otimes _{S}M\overset{%
\vartheta _{M}^{l}}{\simeq }M$).

Before we proceed, we give a number of examples of members in an important
family of semirings which does not include any (non-zero) ring, \emph{i.e.}
every semiring in this family is \emph{proper}.

\begin{definition}
A semiring $(S,+_{S},\cdot _{S},\mathbf{0},\mathbf{1})$ is said to be \emph{%
additively idempotent} iff $\mathbf{1}+\mathbf{1}=\mathbf{1},$ or
equivalently iff $a+a=a$ for every $a\in S.$
\end{definition}

\begin{remark}
\label{add-id}If $(S,+_{S},\cdot _{S},\mathbf{0},\mathbf{1})$ is an
additively idempotent semiring, then 
\begin{equation*}
\mathbf{n}:=\underset{n\text{ times}}{\underbrace{\mathbf{1}+_{S}\cdots +_{S}%
\mathbf{1}}}=\mathbf{1}=\mathbf{1}\cdot _{S}\mathbf{1}\neq \mathbf{0}.
\end{equation*}
\end{remark}

\begin{ex}
Every distributive complete lattice $\mathcal{L}=(L,\vee ,\wedge ,\mathbf{0},%
\mathbf{1})$ is an additively idempotent semiring: $\mathbf{1}+_{\mathcal{L}}%
\mathbf{1=1}\vee \mathbf{1}=\mathbf{1}.$ In particular, for any ring $R,$
the lattice $\mathcal{L}(R):=(\mathrm{Ideal}(R),+,\cdot ,R,0)$ of
(two-sided) ideals of $R$ is distributive and complete whence an additively
idempotent semiring: $\mathbf{1}+_{\mathcal{L}}\mathbf{1}=R+R=R=\mathbf{1}.$
Notice that $\mathcal{L}(R)$ has no non-zero zerodivisors if and only if $%
\{0_{R}\}$ is a prime ideal (\emph{i.e. }$R$ is a prime ring).
\end{ex}

\begin{ex}
$\mathbb{B}=\{0,1\}$ is an additively idempotent semiring with addition: $%
0+0=0$ and $1+1=1$ (called the \emph{Boolean semifield}). Notice that $%
\mathbb{B}\ncong \mathbb{Z}/2\mathbb{Z}.$ The semiring $\mathbb{B}$ has many
applications in automata theory and in switching theory where it is called
the \emph{switching algebra }(\cite[p. 7]{Gol1999}).
\end{ex}

\begin{proposition}
\label{sem-prop}Let $S$ be a semiring.

\begin{enumerate}
\item $\mathbb{S}_{S}$ is complete and cocomplete. In particular, it has
equalizers \emph{(}kernels\emph{)} and coequalizers \emph{(}cokernels\emph{)}%
.

\item $S_{S}$ is a regular generator in $\mathbb{S}_{S}$ \emph{(}in the
sense of \emph{\cite[p. 199]{BW2005})}.

\item If $S$ is commutative, then $(\mathbb{S}_{S},\otimes _{S},S;\mathbf{%
\tau })$ is a symmetric monoidal category with symmetric braiding%
\begin{equation*}
\mathbf{\tau }_{(M,N)}:M\otimes _{S}N\simeq N\otimes _{S}M,\text{ }m\otimes
_{S}n\mapsto n\otimes _{S}m.
\end{equation*}
\end{enumerate}
\end{proposition}

The proofs of the following lemmata are similar to those for modules over a
ring (\emph{e.g.} \cite[12.9, 25.5 (2)]{Wis1991}) by applying a relaxed
version of the Short Five Lemma for semimodules over semirings \cite[Lemma
1.22]{Abu-a}.$\blacksquare $

\begin{definition}
Let $M$ and $N$ be $S$-semimodules. We call an $S$-linear map $%
f:M\longrightarrow N:$

$i$\emph{-uniform} (\emph{image-uniform}) iff 
\begin{equation*}
f(M)=\overline{f(M)}:=\{n\in N\mid n+f(m)=f(m^{\prime })\text{ for some }%
m,m^{\prime }\in M\};
\end{equation*}

$k$\emph{-uniform} (\emph{kernel-uniform}) iff for all $m,m^{\prime }\in M$
we have 
\begin{equation}
f(m)=f(m^{\prime })\Rightarrow m+k=m^{\prime }+k^{\prime }\text{ for some }%
k,k^{\prime }\in \mathrm{Ker}(f);  \label{k-uniform}
\end{equation}

\emph{uniform} iff $f$ is $i$-uniform and $k$-uniform.

We call $L\leq _{S}M$ a \emph{uniform subsemimodule} iff the embedding $L%
\overset{\iota _{L}}{\hookrightarrow }M$ is ($i$-)uniform, or equivalently
iff $L\leq _{S}M$ is subtractive. If $\equiv $ is an $S$-congruence on $M$ 
\cite{Gol1999}, then we call $M/\equiv $ a \emph{uniform quotient} iff the
projection $\pi _{\equiv }:M\longrightarrow M/\equiv $ is ($k$-)uniform.
\end{definition}

\begin{definition}
We say that an $S$-semimodule $X$ (\emph{uniformly}) \emph{generates} $M_{S}$
iff there exists an index set $\Lambda $ and a (uniform) surjective $S$%
-linear map $X^{(\Lambda )}\overset{\pi }{\longrightarrow }M\longrightarrow
0.$ With $\mathrm{Gen}(X)$ we denote the class of $S$-semimodules generated
by $X_{S}.$
\end{definition}

\begin{definition}
We say that $X_{S}$ is

\emph{uniformly\emph{\ }}(\emph{\emph{finitely}})\emph{\ generated} iff
there exists a (finite) index set $\Lambda $ and a uniform surjective $S$%
-linear map $S^{(\Lambda )}\longrightarrow X\longrightarrow 0;$

\emph{finitely presented} iff $\mathrm{Hom}_{S}(X,-):S_{S}\longrightarrow 
\mathbf{AbMonoid}$ preserves directed colimits (\emph{i.e.} $X\in \mathbb{S}%
_{S}$ is a finitely presentable object in the sense of \cite{AP1994});

\emph{uniformly finitely presented} iff $X$ is uniformly finitely generated
and for any exact sequence of $S$-semimodules%
\begin{equation*}
0\longrightarrow K\overset{f}{\longrightarrow }S^{n}\overset{g}{%
\longrightarrow }X\longrightarrow 0,
\end{equation*}%
the $S$-semimodule $K$ ($\simeq \mathrm{Ker}(g)$) is finitely generated.
\end{definition}

\begin{punto}
\label{sgm-M}Let $M$ be a right $S$-semimodule. With $\sigma \lbrack M_{S}]$
($\sigma _{u}[M_{S}]$) we denote the closure of $\mathrm{Gen}(M_{S})$ under (%
\emph{uniform}) $S$-subsemimodules, \emph{i.e.} the smallest full
subcategory of $\mathbb{S}_{S}$ which contains $M_{S}$ and is closed under
direct sums, homomorphic images and (uniform) $S$-subsemimodules. We say
that $M_{S}$ is a (uniformly) \emph{subgenerator} for $\sigma \lbrack M_{S}]$
($\sigma _{u}[M_{S}]$). Notice that $\mathrm{Gen}(M_{S})\subseteq \sigma
_{u}[M_{S}]\subseteq \sigma \lbrack M_{S}].$
\end{punto}

\begin{lemma}
\label{fp-iso}Let $M_{S}$ be a right $S$-semimodule, $\{L_{\lambda
}\}_{\Lambda }$ a class of left $S$-semimodules and consider the canonical
map%
\begin{equation}
\varphi _{M}:M\otimes _{S}\prod\limits_{\lambda \in \Lambda }L_{\lambda
}\longrightarrow \prod\limits_{\lambda \in \Lambda }(M\otimes _{S}L_{\lambda
}),\text{ }m\otimes _{S}(l_{\lambda })_{\Lambda }\mapsto (m\otimes
_{S}l_{\lambda })_{\Lambda }.
\end{equation}

\begin{enumerate}
\item $M_{S}$ is finitely generated if and only if $\varphi _{M}$ is
surjective.

\item If $M_{S}$ is uniformly finitely presented, then $\varphi _{M}$ is an
isomorphism.
\end{enumerate}
\end{lemma}

\begin{definition}
(\cite{Kat2004}, \cite{Abu-a}) We call a right $S$-semimodule $M:$

\emph{flat} iff $M\otimes _{A}-$ is left exact, \emph{i.e.} it preserves
finite limits, equivalently $M\simeq \lim\limits_{\longrightarrow
}F_{\lambda },$ a filtered limit of finitely generated free right $S$%
-semimodules;

\emph{uniformly flat} iff $M\otimes _{A}-:$ $_{A}\mathbb{S}\longrightarrow 
\mathbf{AbMonoid}$ preserves \emph{uniform} subobjects;

\emph{mono-flat} iff $M\otimes _{A}-:$ $_{A}\mathbb{S}\longrightarrow 
\mathbf{AbMonoid}$ preserves monomorphisms (injective $S$-linear maps);

$u$\emph{-flat} iff $M\otimes _{A}-:$ $_{A}\mathbb{S}\longrightarrow \mathbf{%
AbMonoid}$ sends (uniform) monomorphisms to (uniform) monomorphisms;

\emph{projective} iff $M$ is a retract of a free $S$-semimodule, or
equivalently, iff $M$ has a \emph{dual basis}.
\end{definition}

\begin{lemma}
\label{fp-flat}Let $S,T$ be semirings, $L$ a right $S$-semimodule and $K$ a $%
(T,S)$-bisemimodule. Let $Q_{T}$ be a right $T$-semimodule and consider the
canonical morphism%
\begin{equation*}
\upsilon _{(Q,L,K)}:Q\otimes _{T}\mathrm{Hom}_{S}(L,K)\longrightarrow 
\mathrm{Hom}_{S}(L,Q\otimes _{T}K),\text{ }q\otimes _{T}h\mapsto \lbrack
l\mapsto q\otimes h(l)].
\end{equation*}

\begin{enumerate}
\item If $Q_{T}$ is mono-flat and $L_{S}$ is finitely generated, then $%
\upsilon _{(Q,L,K)}$ is injective.

\item If $Q_{T}$ is uniformly flat and $L_{S}$ is uniformly finitely
presented, then $\upsilon _{(Q,L,K)}$ is surjective.

\item If $Q_{T}$ is flat and $L_{S}$ is uniformly finitely presented, then $%
\upsilon _{(Q,L,K)}$ is an isomorphism.
\end{enumerate}
\end{lemma}

\subsection*{Semicoalgebras and Semicomodules}

\qquad Throughout, $S$ is a commutative semiring with $1_{S}\neq 0_{S}.$

\begin{punto}
An $S$\emph{-semialgebra}\ is a triple $(A,\mu _{A},\eta _{A})$ where $A$ is
an $S$-semimodule and $\mu _{A}:A\otimes _{S}A\longrightarrow A,$ $\eta
_{A}:S\longrightarrow A$ are $S$-linear maps such that the following
diagrams are commutative%
\begin{equation*}
\begin{array}{ccc}
\xymatrix{A \otimes_S A \otimes_S A \ar[rr]^{\mu_A \otimes_S A} \ar[dd]_{A
\otimes_S \mu_A} & & A \otimes_S A \ar[dd]^{\mu_A} \\ & & \\ A \otimes_S A
\ar[rr]_{\mu_A} & & A} &  & \xymatrix{S \otimes_S A \ar[ddrr]_{\vartheta_A
^l} \ar[rr]^{\eta_A \otimes_S A } & & A \otimes_S A \ar[dd]_{\mu_A} & & A
\otimes_S S \ar[ll]_{A\otimes_S \eta_A} \ar[ddll]^{\vartheta_A ^r} \\ & & &
& & \\ & & A & & }%
\end{array}%
\end{equation*}%
We call $\mu _{A}$ the \emph{multiplication }and $\eta _{A}$ the \emph{unity}
of $A.$ Let $A$ and $B$ be $S$-semialgebras. We call an $S$-linear map $%
f:A\longrightarrow B$ an $S$\emph{-semialgebra morphism }iff the following
diagrams are commutative%
\begin{equation*}
\begin{array}{ccc}
\xymatrix{A \otimes_S A \ar^(.45){f \otimes_S f}[rr] \ar_(.45){\mu_A}[d] & &
B \otimes_S B \ar^(.45){\mu_B}[d] \\ A \ar_(.45){f}[rr] & & B} &  & %
\xymatrix{A \ar[rr]^{f} & & B \\ & S \ar[ur]_{\eta_B} \ar[ul]^{\eta_A} &}%
\end{array}%
\end{equation*}%
The set of morphisms of $S$-semialgebras form $A$ to $B$ is denoted by $%
\mathrm{SAlg}_{S}(A,B).$ The category of $S$-semialgebras will be denoted by 
$\mathbf{SAlg}_{S}.$
\end{punto}

Semicoalgebras are dual to semialgebras and are defined by reversing the
arrows in the diagrams mentioned above.

\begin{punto}
An\emph{\ }$S$\emph{-semicoalgebra} is a triple $(C,\Delta _{C},\varepsilon
_{C})$ in which $C$ is an $S$-semimodule and $\Delta _{C}:C\longrightarrow
C\otimes _{S}C,$ $\varepsilon _{C}:C\longrightarrow S$ are $S$-linear maps
such that the following diagrams are commutative%
\begin{equation*}
\begin{array}{ccc}
\xymatrix{C \ar^(.45){\Delta_C}[rr] \ar_(.45){\Delta_C}[d] & & C \otimes_S C
\ar^(.45){C \otimes_S \Delta_C}[d]\\ C \otimes_S C \ar_(.45){\Delta_C
\otimes_S C}[rr] & & C \otimes_S C \otimes_S C } &  & \xymatrix{ & C
\ar^(.45){\Delta _C}[d] & \\ {S} \otimes_{S} C \ar[ur]^(.45){\vartheta _C
^l} & C \otimes_{S} C \ar^(.45){\varepsilon _{C} \otimes_S C}[l] \ar_(.45){C
\otimes_S \varepsilon _{C}}[r] & C \otimes_{S} {S} \ar[ul]_(.45){\vartheta_C
^r} }%
\end{array}%
\end{equation*}%
We call $\Delta _{C}$ the \emph{comultiplication} and $\varepsilon _{C}$ the 
\emph{counity} of $C.$ For $S$-semicoalgebras $C$ and $D,$ we call an $S$%
-linear map $f:D\longrightarrow C$\ an $S$\emph{-semicoalgebra morphism} iff
the following diagrams are commutative%
\begin{equation*}
\begin{array}{ccc}
\xymatrix{D \ar^(.45){f}[rr] \ar_(.45){\Delta_D}[d] & & C
\ar^(.45){\Delta_{C}}[d] \\ D \otimes_{S} D \ar_(.45){f \otimes_S f}[rr] & &
C \otimes_{S} C & &} &  & \xymatrix{D \ar[dr]_{\varepsilon_D} \ar[rr]^{f} &
& C \ar[dl]^{\varepsilon_C} \\ & S &}%
\end{array}%
\end{equation*}%
The set of $S$-semicoalgebra morphisms from $D$ to $C$ is denoted by $%
\mathrm{SCoalg}_{S}(D,C).$ The category of $S$-semicoalgebras is denoted by $%
\mathbf{SCoalg}_{S}.$
\end{punto}

\begin{notation}
Let $(C,\Delta ,\varepsilon )$ be an $S$-semicoalgebra. We use
Sweedler-Heyneman's $\sum $-notation, and write for $c\in C:$%
\begin{eqnarray*}
\Delta (c) &=&\sum c_{1}\otimes _{S}c_{2}\in C\otimes _{S}C; \\
\sum c_{11}\otimes _{S}c_{12}\otimes _{S}c_{2} &=&\sum c_{1}\otimes
_{S}c_{2}\otimes _{S}c_{3}=\sum c_{1}\otimes _{S}c_{21}\otimes _{S}c_{22}.
\end{eqnarray*}
\end{notation}

\begin{punto}
Notice that an $S$-semialgebra $(A,\mu ,\eta )$ is \emph{commutative} iff $%
\mu _{A}\circ \mathbf{\tau }_{(A,A)}=\mu _{A};$ with $_{\mathrm{c}}\mathbf{%
SAlg}_{S}\hookrightarrow \mathbf{SAlg}_{S},$ we denote the category of
commutative $S$-semialgebras. Dually, an $S$-semicoalgebra $(C,\Delta
,\varepsilon )$ is said to be \emph{cocommutative} iff $\mathbf{\tau }%
_{(C,C)}\circ \Delta =\Delta ,$ \emph{i.e. }$\sum c_{1}\otimes
_{S}c_{2}=\sum c_{2}\otimes _{S}c_{1}$ for all $c\in C.$ With $_{\mathrm{coc}%
}\mathbf{SCoalg}_{S}\hookrightarrow \mathbf{SCoalg}_{S},$ we denote the 
\emph{full} subcategory of cocommutative $S$-semicoalgebras.
\end{punto}

\begin{ex}
Let $M$ be an $S$-semimodule. We have an $S$-semicoalgebra structure on $%
C=(S\oplus M,\Delta ,\varepsilon ),$ where%
\begin{eqnarray*}
\Delta &:&(s,m)\mapsto (s,0)\otimes _{S}(1,0)+(1,0)\otimes
_{S}(0,m)+(0,m)\otimes _{S}(1,0); \\
\varepsilon &:&(s,m)\mapsto s.
\end{eqnarray*}%
Notice that there are many properties $\mathbb{P}$ such that $_{S}C$ has
Property $\mathbb{P}$ if (and only if) $_{S}M$ has Property $\mathbb{P},$ 
\emph{e.g.} being flat, (finitely) projective, finitely generated \cite[%
Example 10 (1)]{Wisch1975}.
\end{ex}

\begin{ex}
Let $X$ be any set. We have an $S$-semicoalgebra $(S[X],\Delta ,\varepsilon
),$ where $S[X]$ is the free $S$-semimodule with basis $X$ and $\Delta
,\varepsilon $ are defined by extending the following assignments linearly%
\begin{equation*}
\Delta :S[X]\mapsto S[X]\otimes _{S}S[X],\text{ }x\mapsto x\otimes _{S}x%
\text{ and }\varepsilon :S[X]\mapsto S,\text{ }x\mapsto 1_{S}.
\end{equation*}
\end{ex}

\subsection*{Semicomodules}

\qquad Dual to semimodules of semialgebras are semicomodules of
semicoalgebras:

\begin{punto}
Let $(C,\Delta ,\varepsilon )$ be an $S$-semicoalgebra. A right $C$\emph{%
-semicomodule} is an $S$-semimodule $M$ associated with an $S$-linear map
(called $C$\emph{-coaction})%
\begin{equation*}
\rho ^{M}:M\longrightarrow M\otimes _{S}C,\text{ }m\mapsto \sum
m_{<0>}\otimes _{S}m_{<1>},
\end{equation*}%
such that the following diagrams are commutative%
\begin{equation*}
\begin{array}{ccc}
\xymatrix{M \ar^(.4){\rho ^M}[rr] \ar_(.45){\rho ^M}[d] & & M \otimes_S {C}
\ar^(.45){M \otimes_S \Delta}[d] \\ M \otimes_S {C} \ar_(.4){\rho ^M
\otimes_S {C}}[rr] & & M \otimes_S C \otimes_S {C}} &  & \xymatrix{M
\ar[rr]^{\rho^M} & & M \otimes_S C \ar[dl]^{M \otimes_S \varepsilon} \\ & M
\otimes_S S \ar[ul]^{\vartheta_M^r} & }%
\end{array}%
\end{equation*}%
Let $M$ and $N${\normalsize \ }be right $C$-semicomodules. We call an $S$%
-linear map $f:M\longrightarrow N$ a $C$\emph{-semicomodule morphism}%
{\normalsize \ }(or $C$\emph{-colinear}) iff the following diagram is
commutative%
\begin{equation*}
\xymatrix{M \ar[rr]^{f} \ar[d]_{\rho ^M} & & N \ar[d]^{\rho ^N}\\ M
\otimes_S {C} \ar[rr]_{f \otimes_S {C}} & & N \otimes_S {C} }
\end{equation*}%
The set of $C$-colinear maps from $M$ to $N$ is denoted by $\mathrm{Hom}%
^{C}(M,N).$ The category of right $C$-semicomodules and $C$-colinear maps is
denoted by $\mathbb{S}^{C}.$ For a right $C$-semicomodule $M,$ we call $%
L\leq _{A}M$ a $C$\emph{-subsemicomodule} iff $(L,\rho ^{L})\in \mathbb{S}%
^{C}$ and the embedding $L\overset{\iota _{L}}{\hookrightarrow }M$ is $C$%
-colinear. Symmetrically, we define the category $^{C}\mathbb{S}$ of left $C$%
-semicomodules. For two left $C$-semicomodules $M$ and $N,$ we denote by $%
^{C}\mathrm{Hom}(M,N)$ the set of $C$-colinear maps from $M$ to $N.$
\end{punto}

\begin{punto}
Let $(M,\rho ^{(M;C)})$ be a right $C$-semicomodule, $(M,\rho ^{(M;D)}$ a
left $D$-semicomodule and consider the left $D$-semicomodule $(M\otimes
_{S}C,\rho ^{(M;D)}\otimes _{S}C)$ (the right $C$-semicomodule $(D\otimes
_{S}M,D\otimes _{S}\rho ^{(M;C)}$). We call $M$ a $(D,C)$\emph{%
-bisemicomodule} iff $\rho ^{(M;C)}:M\longrightarrow M\otimes _{S}C$ is $D$%
-colinear, or equivalently iff $\rho ^{(M;D)}:M\longrightarrow D\otimes
_{S}M $ is $C$-colinear. For $(D,C)$-bisemicomodules $M$ and $N,$ we call a $%
D$-colinear $C$-colinear map $f:M\longrightarrow N$ a $(D,C)$\emph{%
-bisemicomodule morphism}\textbf{\ }(or $(D,C)$\emph{-bicolinear}). The
category of $(D,C)$-bisemicomodules and $(D,C)$-bicolinear maps is denoted
by $^{D}\mathbb{S}^{C}.$
\end{punto}

\begin{remark}
\label{inj-kounit}Let $(C,\Delta ,\varepsilon )$ be an $S$-semicoalgebra. If 
$(M,\rho ^{M})$ is a right $C$-semicomodule, then $\rho
^{M}:M\longrightarrow M\otimes _{S}C$ is a splitting monomorphism in $%
\mathbb{S}_{A};$ however, $M$ is not necessarily a direct summand of $%
M\otimes _{S}C\mathbf{;}$ see \cite[16.6]{Gol1999}.
\end{remark}

\begin{ex}
Let $M=\dbigoplus\limits_{g\in G}M_{g}$ be a $G$-graded $S$-semimodules,
where $G$ is group. One can consider $M$ as an $S[G]$-semimodule with%
\begin{equation*}
\rho ^{M}:M\longrightarrow M\otimes _{S}S[G],\text{ }\sum_{g\in
G}m_{g}\mapsto \sum_{g\in G}m_{g}\otimes _{S}g.
\end{equation*}%
Conversely, if $M$ is an $S[G]$-semicomodule, then $M$ is a $G$-graded $S$%
-semimodule in the canonical way. In fact, we have an isomorphism of
categories $\mathbf{gr}_{G}(\mathbb{S}_{S})\simeq \mathbb{S}^{S[G]},$ where $%
\mathbf{gr}_{G}(\mathbb{S}_{S})$ is the category of $G$-graded $S$%
-semimodules.
\end{ex}

\begin{punto}
We have an isomorphism of categories $\mathbf{SAlg}_{S}\simeq \mathbf{Monoid}%
(\mathbb{S}_{S}).$ An $S$-semialgebra $A$ is essentially a \emph{monoid} in $%
\mathbb{S}_{S}$ and so it induces two \emph{monads} $-\otimes _{S}A:\mathbb{S%
}_{S}\longrightarrow \mathbb{S}_{S}$ and $A\otimes _{S}-:\mathbb{S}%
_{S}\longrightarrow \mathbb{S}_{S}.$ Moreover, we have isomorphisms of
categories%
\begin{equation*}
\mathbb{S}_{A}\simeq (\mathbb{S}_{S})_{-\otimes _{S}A}\text{ and }_{A}%
\mathbb{S}\simeq (\mathbb{S}_{S})_{A\otimes _{S}-}.
\end{equation*}%
We have an isomorphism of categories $\mathbf{SCAlg}_{S}\simeq \mathbf{%
Comonoid}(\mathbb{S}_{S}).$ An $S$-semicoalgebra $C$ is essentially a \emph{%
comonoid} in $\mathbb{S}_{S}$ and so it induces two \emph{comonads} $%
-\otimes _{S}C:\mathbb{S}_{S}\longrightarrow \mathbb{S}_{S}$ and $C\otimes
_{S}-:\mathbb{S}_{S}\longrightarrow \mathbb{S}_{S}.$ Moreover, we have
isomorphisms of categories%
\begin{equation*}
\mathbb{S}^{C}\simeq \mathbb{S}_{S}^{-\otimes _{S}C}\text{ and }^{C}\mathbb{S%
}\simeq \mathbb{S}_{S}^{C\otimes _{S}-}.
\end{equation*}
\end{punto}

\begin{punto}
Let $C$ be an $S$-semicoalgebra. For every $S$-semialgebras $A,$ there is a
canonical structure of an $S$-semialgebra on $\mathrm{Hom}_{S}(C,A)$ with
multiplication given by the \emph{convolution product}%
\begin{equation}
(f\ast g)(c)=\sum f(c_{1})g(c_{2})\text{ for all }f,g\in \mathrm{Hom}%
_{S}(C,A)\text{ and }c\in C.  \label{conv}
\end{equation}%
In particular, $C^{\ast }:=\mathrm{Hom}_{S}(C,S)$ is an $S$-semialgebra and $%
C$ is a $(C^{\ast },C^{\ast })$-bisemimodule with left and right actions
given by%
\begin{equation*}
f\rightharpoonup c:=\sum c_{1}f(c_{2})\text{ and }c\leftharpoonup g:=\sum
f(c_{1})c_{2}\text{ for all }f,g\in C^{\ast }\text{ and }c\in C.
\end{equation*}
\end{punto}

\begin{punto}
Let $C$ be an $S$-semicoalgebra. We say that $_{S}C$ is an $\alpha $\emph{%
-semimodule} (or $_{S}C$ \emph{satisfies the }$\alpha $\emph{-condition})
iff for every $M_{S},$ the canonical map%
\begin{equation*}
\alpha _{M}^{C}:M\otimes _{S}C\longrightarrow \mathrm{Hom}_{S}(C^{\ast },M),%
\text{ }m\otimes _{S}c\mapsto \lbrack f\mapsto mf(c)]
\end{equation*}%
is injective and uniform. Clearly, every right $C$-semicomodule $M$ is a
left $^{\ast }C$-semimodule with%
\begin{equation*}
f\rightharpoonup m:=\sum m_{<0>}f(m_{<1>})\text{ for all }f\in C^{\ast }%
\text{ and }m\in M.
\end{equation*}%
If $_{S}C$ is an $\alpha $-semimodule, then for every $_{^{\ast }C}M$ with
induced map $\widetilde{\rho }_{M}:M\longrightarrow \mathrm{Hom}_{S}(C^{\ast
},M),$ we define the $C$\emph{-rational subsemimodule} of $M$ as $\mathrm{Rat%
}^{C}(_{^{\ast }C}M):=\widetilde{\rho }_{M}^{-1}(\alpha _{M}^{C}(M\otimes
_{S}C)).$ In this case, we have by \cite[Theorem 3.16]{Abu-c} an isomorphism
of categories 
\begin{equation}
\mathbb{S}^{C}\simeq \mathrm{Rat}^{C}(_{^{\ast }C}\mathbb{S}).
\label{rat-co}
\end{equation}
\end{punto}

\section{Bisemialgebras}

\begin{punto}
With an $S$\emph{-bisemialgebra}, we mean a datum $(B,\mu ,\eta ,\Delta
,\varepsilon ),$ where $(B,\mu ,\eta )$ is an $S$-semialgebra and $(B,\Delta
,\varepsilon )$ is an $S$-semicoalgebra such that $\Delta :B\longrightarrow
B\otimes _{S}B$ and $\varepsilon :B\longrightarrow S$ are morphisms of $S$%
-semialgebras, or equivalently $\mu :B\otimes _{S}B\longrightarrow B$ and $%
\eta :S\longrightarrow B$ are morphisms of $S$-semicoalgebras; notice that $%
B\otimes _{S}B$ can be given a structure of a $S$-semialgebra ($S$%
-semicoalgebra) in a canonical way using the twisting map $\mathbf{\tau }%
_{(B,B)}.$ A \emph{morphism of }$S$\emph{-bisemialgebras }$%
f:B\longrightarrow B^{\prime }$ is an $S$-linear map which is simultaneously
a morphism of $S$-semialgebras and a morphism of $S$-semicoalgebras. The
category of $S$-bisemialgebras is denoted by $\mathbf{SBiAlg}_{S}.$
\end{punto}

\begin{notation}
Given an $S$-bisemialgebra $B,$ we write $B^{a}$ when we handle $B$ as an $S$%
-semialgebra and $B^{c}$ when we consider $B$ as an $S$-semicoalgebra.
\end{notation}

\begin{ex}
$S$ is an $S$-bisemialgebra with%
\begin{eqnarray*}
\Delta _{S} &:&S\longmapsto S\otimes _{S}S,\otimes _{S},\text{ }s\mapsto
s\otimes _{S}1_{S}=1_{S}\otimes _{S}s; \\
\varepsilon _{S} &:&S\longrightarrow S,\text{ }s\mapsto s.
\end{eqnarray*}
\end{ex}

\begin{ex}
If $B$ is an $S$-bisemialgebra, then $S$ is a $(B,B)$-bisemicomodule with%
\begin{equation*}
^{S}\rho :S\longrightarrow B\otimes _{S}S,\text{ }s\longmapsto 1_{B}\otimes
_{S}s\text{ and }\rho ^{S}:S\longrightarrow S\otimes _{S}B,\text{ }%
s\longmapsto s\otimes _{S}1_{B}.
\end{equation*}
\end{ex}

\begin{ex}
Let $A$ be an $S$-semialgebra and consider $B=S\oplus A$ as an $S$%
-semialgebra with point wise multiplication and unity $(1_{S},1_{A}).$ It is
obvious that $B$ has a structure of an $S$-bisemialgebra with%
\begin{eqnarray*}
\Delta &:&B\longrightarrow B\otimes _{S}B,\text{ }(s,a)\mapsto (s,0)\otimes
_{S}(1_{S},0)+(1_{S},0)\otimes _{S}(0,a)+(0,a)\otimes _{S}(1_{S},0); \\
\varepsilon &:&B\longrightarrow S,\text{ }(s,a)\mapsto s.
\end{eqnarray*}
\end{ex}

\begin{ex}
\label{S[M]}Let $(G,\ast ,e)$ be a monoid and consider the free $S$%
-semimodule $S[G]$ as an $S$-semialgebra with multiplication induced by $%
\ast $ and unity $1=1_{S}e.$ One can easily see that $S[G]$ has two $S$%
-semicoalgebra structures which are compatible with the $S$-semialgebra
structure yielding two $S$-bisemialgebra structures $(S[G],\ast ,1,\Delta
,\varepsilon )$ and $(S[G],\ast ,1,\Delta ,\varepsilon )$ where the
comultiplications and the counities are obtained by extending the following
assignments as $S$-semialgebra morphisms%
\begin{eqnarray*}
\Delta &:&S[G]\longrightarrow S[G]\otimes _{S}S[G],\text{ }g\mapsto g\otimes
_{S}g\text{ and }\varepsilon :S[G]\longrightarrow S,\text{ }g\mapsto 1_{S};
\\
\widetilde{\Delta } &:&S[G]\longrightarrow S[G]\otimes _{S}S[G],\text{ }%
g\mapsto g\otimes _{S}1+1\otimes _{S}g\text{ and }\widetilde{\varepsilon }%
:S[G]\longrightarrow S,\text{ }g\mapsto \delta _{e,g}.
\end{eqnarray*}
\end{ex}

\begin{ex}
The previous example applies in particular to the polynomial $S$-semialgebra 
$S[x]$ since we have an isomorphism of monoids $M:=(\{1,X,\cdots
,X^{n},\cdots \},\cdot )\simeq (\mathbb{N}_{0},+),$ whence $S[x]$ has two $S$%
-bisemialgebra structures $(S[x],\cdot ,1,\Delta _{1},\varepsilon _{1})$ and 
$(S[x],\cdot ,1,\Delta _{2},\varepsilon _{2})$ with%
\begin{eqnarray*}
\Delta &:&S[x]\longrightarrow S[x]\otimes _{S}S[x],\text{ }%
\sum_{i=0}^{n}s_{i}x^{i}\mapsto \sum_{i=0}^{n}s_{i}x^{i}\otimes _{S}x^{i}; \\
\varepsilon &:&S[x]\longrightarrow S,\text{ }\sum_{i=0}^{n}s_{i}x^{i}\mapsto
\sum_{i=0}^{n}s_{i};
\end{eqnarray*}%
and%
\begin{eqnarray*}
\widetilde{\Delta } &:&S[x]\longrightarrow S[x]\otimes _{S}S[x],\text{ }%
\sum_{i=0}^{n}s_{i}x^{i}\mapsto \sum_{i=0}^{n}s_{i}\left( \sum_{j=0}^{i}%
\binom{i}{j}x^{j}\otimes _{S}x^{j-i}\right) ; \\
\widetilde{\varepsilon } &:&S[x]\longrightarrow S,\text{ }%
\sum_{i=0}^{n}s_{i}x^{i}\mapsto s_{0}.
\end{eqnarray*}
\end{ex}

\begin{ex}
\label{Worth}(cf. \cite{Wor2012}) Consider the Boolean semiring $\mathbf{B}%
=\{0,1\}.$ Let $P=\mathbf{B}<x,y\mid xy\neq yx>$ be the $\mathbf{B}$%
-semimodule of formal sums of words formed from the \emph{non-commuting}
letters $x$ and $y.$ In fact, $P$ is a non-commutative $\mathbf{B}$%
-semialgebra with multiplication given by \emph{concatenation} of words (%
\emph{e.g. }$(xyxx)\cdot (yyx)=xyxxyyx$) and unity $[]$ (the empty word). It
can be easily seen that the structure maps of the following $\mathbf{B}$%
-semicoalgebras can be extended as $S$-semialgebra morphisms yielding two
different $\mathbf{B}$-bisemialgebra structures on $P:$

\begin{enumerate}
\item $(P,\Delta ,\varepsilon ),$ where%
\begin{equation*}
\Delta _{1}:P\longrightarrow P\otimes _{\mathbf{B}}P,\text{ }w\mapsto
w\otimes _{\mathbf{B}}w\text{ and }\varepsilon _{1}:P\longrightarrow \mathbf{%
B},\text{ }w\mapsto w(1,1).
\end{equation*}

\item $(P,\widetilde{\Delta },\widetilde{\varepsilon }),$ where%
\begin{eqnarray*}
\widetilde{\Delta } &:&P\longrightarrow P\otimes _{\mathbf{B}}P,\text{ }%
x\mapsto x\otimes _{\mathbf{B}}[]+[]\otimes _{\mathbf{B}}x,\text{ }\Delta
(y)=y\otimes _{\mathbf{B}}[]+[]\otimes _{\mathbf{B}}y; \\
\widetilde{\varepsilon } &:&P\longrightarrow \mathbf{B},\text{ }w\mapsto
w(0,0)\text{ (notice that }[](0,0)=1_{S}\text{ by convention). }
\end{eqnarray*}
\end{enumerate}
\end{ex}

\begin{ex}
\label{Haz}(cf. \cite[3.2.20]{HGK2010}, \cite{Wor2012})\ Let $S$ be a
semiring and consider the free $S$-semimodule $B=S<\mathbb{N}>$ with basis
all \emph{words} on $\mathbb{N}.$ Notice that $B$ is an $S$-semialgebra with
multiplication given by the \emph{concatenation} of words and unity $[],$
the empty word. Moreover, $B$ is an $S$-semicoalgebra with comultiplication
and counity given by%
\begin{equation*}
\Delta :B\longrightarrow B\otimes _{S}B,\text{ }w\mapsto
\sum_{w_{1}w_{2}=w}w_{1}\otimes _{S}w_{2}\text{ and }\varepsilon
:B\longrightarrow S,\text{ }w\mapsto w(0,,\cdots ,0).
\end{equation*}%
However, these $S$-semialgebra and $S$-semicoalgebra structures are -- in
general --\emph{\ not} compatible and so do not yield a structure of an $S$%
-bisemialgebra on $B;$ for example, we have%
\begin{eqnarray*}
\Delta ([2])\cdot \Delta ([3]) &=&([]\otimes _{R}[2]+[2]\otimes _{R}[])\cdot
([]\otimes _{R}[3]+[3]\otimes _{R}[]) \\
&=&[]\otimes _{R}[2,3]+[2]\otimes _{R}[3]+[3]\otimes _{R}[2]+[2,3]\otimes
_{R}[],
\end{eqnarray*}%
while%
\begin{equation*}
\Delta ([2,3])=[]\otimes _{R}[2,3]+[2]\otimes _{R}[3]+[2,3]\otimes _{R}[].
\end{equation*}
\end{ex}

\begin{ex}
\label{E}(cf. \cite[Example 7.9]{Str2007}) Let the commutative semiring $S$
be additively idempotent, $E=\{e_{0},e_{1},e_{2},\cdots ,e_{n},\cdots ,\}$ a
countable set, $B:=S[E]$ the free $S$-semimodule and consider the assignments%
\begin{eqnarray*}
\mu _{E} &:&B\otimes _{S}B\longrightarrow B,\text{ }e_{p}\otimes
_{S}e_{q}\mapsto e_{p+q}; \\
\eta &:&S\longrightarrow B,\text{ }1_{S}\mapsto e_{0}; \\
\Delta &:&B\longrightarrow B\otimes _{S}B,\text{ }e_{n}\mapsto
\sum_{p+q=n}e_{p}\otimes _{S}e_{q}; \\
\varepsilon &:&B\longrightarrow S,\text{ }e_{n}\mapsto \delta _{0,n}.
\end{eqnarray*}%
Extending these assignments linearly such that $\mu $ and $\eta $ are $S$%
-semialgebra morphisms, we obtain a structure of an $S$-bisemialgebra on $B.$
\end{ex}

\begin{ex}
(\cite[15.12]{BW2003})\ Associated to every $S$-semimodule $M$ is the
so-called $S$\emph{-tensor-semialgebra}%
\begin{equation*}
\mathcal{T}(M)=(S\oplus M\oplus (M\otimes _{S}M)\oplus (M\otimes
_{S}M\otimes _{S}M)\oplus \cdots ,\mu ,\eta )
\end{equation*}%
where the multiplication and the unity are given by%
\begin{eqnarray*}
\mu \left( \left( m_{1}\otimes _{S}\cdots \otimes _{S}m_{n}\right) \left(
m_{1}^{\prime }\otimes _{S}\cdots \otimes _{S}m_{t}^{\prime }\right) \right)
&:&=m_{1}\otimes _{S}\cdots \otimes _{S}m_{n}\otimes _{S}m_{1}^{\prime
}\otimes _{S}\cdots \otimes _{S}m_{t}^{\prime }; \\
\eta (s) &:&=(s,0,0,0,\cdots ).
\end{eqnarray*}%
Notice that $M\hookrightarrow \mathcal{T}(M),$ $m\mapsto (0,m,0,0,\cdots ).$
In fact, we have an adjoint pair of functors $(\mathcal{T}(-),\mathcal{U}),$
where $\mathcal{U}:\mathbf{SAlg}_{S}\longrightarrow \mathbb{S}_{S}$ is the
forgetful functor. In other words, $\mathcal{T}(M)$ satisfies the following
universal property: given an $S$-linear map $g:M\longrightarrow A,$ where $A$
is an $S$-semialgebra, there exists a morphism of $S$-semialgebras $%
\widetilde{g}:\mathcal{T}(M)\longrightarrow A$ such that $\widetilde{g}\circ
\iota =g.$ By this universal property, the $S$-linear maps%
\begin{equation*}
g:M\longrightarrow \mathcal{T}(M)\otimes _{S}\mathcal{T}(M),\text{ }m\mapsto
m\otimes _{S}1+1\otimes _{S}m\text{ and }\mathfrak{z}:M\longrightarrow S,%
\text{ }m\mapsto 0,
\end{equation*}%
induce $S$-semialgebra morphisms%
\begin{equation*}
\Delta :\mathcal{T}(M)\longrightarrow \mathcal{T}(M)\otimes _{S}\mathcal{T}%
(M)\text{ and }\varepsilon :\mathcal{T}(M)\longrightarrow S,\text{ }m\mapsto
0.
\end{equation*}%
One can easily check that $(\mathcal{T}(M),\mu ,\eta ,\Delta ,\varepsilon )$
is an $S$-bisemialgebra.
\end{ex}

\begin{punto}
Let $C$ be an $S$-semicoalgebra. With a coideal of $C,$ we mean an $S$%
-subsemimodule $K\leq _{S}C$ such that $K=\mathrm{Ker}(f)$ for some \emph{%
uniform} surjective morphism of $S$-semicoalgebras $f:C\longrightarrow
C^{\prime }.$ For characterizations of coideals of semicoalgebra over
semirings, see \cite[Proposition 2.16]{Abu-c}.
\end{punto}

\begin{definition}
A \emph{bi-ideal} $I$ of an $S$-bisemialgebra $B$ is an ideal of $B^{a}$
which is also a coideal of $B^{c}.$
\end{definition}

\begin{ex}
Let $B$ be an $S$-bisemialgebra such that $\varepsilon _{B}$ is uniform.
Notice that $\varepsilon _{B}:B\longrightarrow S$ is a surjective morphism
of $S$-bisemialgebras, whence $\mathrm{Ker}(\varepsilon _{B})$ is a bi-ideal.
\end{ex}

\begin{lemma}
Let $\gamma :B\longrightarrow B^{\prime }$ be a morphism of $S$%
-bisemialgebras.

\begin{enumerate}
\item If $\gamma $ is surjective and uniform, then $\mathrm{Ker}(\gamma )$
is a bi-ideal.

\item $\gamma (B)$ is an $S$-subbisemialgebra of $B^{\prime }.$

\item For any bi-ideal $I\subseteq \mathrm{Ker}(\gamma ),$ there is a
commutative diagram of $S$-bisemialgebras%
\begin{equation*}
\xymatrix{B \ar[rr]^{\gamma} \ar[dr]_{\pi_I} & & B' \\ & B/I \ar[ur]_{f} & }
\end{equation*}%
where $\pi :B\longrightarrow B/I$ is the canonical projection.
\end{enumerate}
\end{lemma}

\subsection*{Hopf Semimodules}

\qquad In what follows, let $(B,\mu ,\eta ,\Delta ,\varepsilon )$ be an $S$%
-bisemialgebra.

\begin{punto}
If $M,N\in \mathbb{S}_{B},$ then $M\otimes _{S}N$ has a \emph{trivial}
structure of a right $B$-semimodule $(M\otimes _{S}N,-\otimes _{S}\rho ^{N})$
and another structure of a right $B$-semimodule $(M\otimes _{S}^{a}N,\rho
_{M\otimes _{S}^{a}N})$ where $\rho _{M\otimes _{S}^{a}N}$ is the
composition of the following maps%
\begin{equation*}
(M\otimes _{S}^{a}N)\otimes _{S}B\overset{-\otimes _{S}\Delta }{%
\longrightarrow }(M\otimes _{S}^{a}N)\otimes _{S}(B\otimes _{S}B)\overset{%
-\otimes _{S}\mathbf{\tau }_{(N,B)}\otimes _{S}-}{\simeq }(M\otimes
_{S}B)\otimes _{S}(N\otimes _{S}B)\overset{\rho _{M}\otimes _{S}\rho _{N}}{%
\longrightarrow }M\otimes _{S}^{a}N.
\end{equation*}%
On the other hand, if $M,N\in \mathbb{S}^{B},$ then $M\otimes _{S}N$ has a 
\emph{trivial} structure of a right $B$-semicomodule $(M\otimes
_{S}N,-\otimes _{S}\rho ^{N})$ and another structure of a right $B$%
-semicomodule $(M\otimes _{S}^{c}N,\rho ^{M\otimes _{S}^{c}N})$ where $\rho
^{M\otimes _{S}^{c}N}$ is the composition of the following maps%
\begin{equation*}
M\otimes _{S}^{c}N\overset{\rho ^{M}\otimes _{S}\rho ^{N}}{\longrightarrow }%
(M\otimes _{S}B)\otimes _{S}(N\otimes _{S}B)\overset{-\otimes _{S}\mathbf{%
\tau }_{(B,N)}\otimes -}{\longrightarrow }(M\otimes _{S}^{c}N)\otimes
_{S}(B\otimes _{S}B)\overset{-\otimes _{S}\mu }{\longrightarrow }(M\otimes
_{S}^{c}N)\otimes _{S}B.
\end{equation*}
\end{punto}

\begin{punto}
A \emph{right-right Hopf semimodule} over $B$ is a triple $(M,\rho _{M},\rho
^{M})$ such that $(M,\rho _{M})$ is a right $B$-semimodule, $(M,\rho ^{M})$
is a right $B$-semicomodule and $\rho _{M}:M\otimes _{S}^{c}B\longrightarrow
M$ is a morphism of right $B$-semicomodules, or equivalently $\rho
^{M}:M\longrightarrow M\otimes _{S}^{a}B$ is a morphism of right $B$%
-semimodules, \emph{i.e.}%
\begin{equation*}
\sum (mb)_{<0>}\otimes _{S}(mb)_{<1>}=\sum m_{<0>}b_{1}\otimes
_{S}m_{<1>}b_{2}\text{ for all }m\in M,b\in B.
\end{equation*}%
The category of right-right Hopf $B$-semimodules with arrows being the $B$%
-linear $B$-colinear maps is denoted by $\mathbb{S}_{B}^{B},$ \emph{i.e.} $%
\mathrm{Hom}_{B}^{B}(M,N)=\mathrm{Hom}_{B}(M,N)\cap \mathrm{Hom}^{B}(M,N)$
for all $M,N\in \mathbb{S}_{B}^{B}$. Symmetrically, one can define the
category $^{B}\mathbb{S}_{B}$ of \emph{right-left Hopf semimodules}, the
category $_{B}^{B}\mathbb{S}$ of \emph{left-left Hopf }$B$\emph{-semimodules}
and the category $_{B}\mathbb{S}^{B}$ of \emph{left-right Hopf }$B$\emph{%
-semimodules}.
\end{punto}

\begin{remarks}
Let $B$ be an $S$-bisemialgebra.

\begin{enumerate}
\item $B\otimes _{S}^{a}B$ and $B\otimes _{S}^{c}$ are subgenerators in $%
\mathbb{S}_{B}^{B}$ (cf. \cite[14.5]{BW2003}).

\item Let $M,N\in \mathbb{S}_{B}^{B}.$ Since $\mathbb{S}_{S}$ has
equalizers, we have 
\begin{equation*}
\mathrm{Equal}(\varphi ,\psi )=\mathrm{Hom}_{B}^{B}(M,N)=\mathrm{Equal}%
(\varkappa ,\omega )
\end{equation*}%
where%
\begin{equation*}
\varphi (f)=\rho ^{N}\circ f,\text{ }\psi (f)=(f\otimes _{S}B)\circ \rho
^{M},\text{ }\varkappa (g)=\rho _{N}\circ (g\otimes _{S}B)\text{ and }\omega
(g)=g\circ \rho _{M}.
\end{equation*}%
\begin{equation*}
\begin{tabular}{l}
$\xymatrix{ \mathrm{Hom}_{B}^{B}(M,N) \ar[r] & \mathrm{Hom}_{B}(M,N)
\ar@<1ex>^{\varphi}[rr] \ar_{\psi}[rr] & & \mathrm{Hom}_{B}(M,N\otimes
_{S}^{a}B)}$ \\ 
\\ 
$\xymatrix{ \mathrm{Hom}_{B}^{B}(M,N) \ar[r] & \mathrm{Hom}_{B}(M,N)
\ar@<1ex>^{\varkappa}[rr] \ar_{\omega}[rr] & & \mathrm{Hom}_{B}(M,N\otimes
_{S}^{c}B)}$%
\end{tabular}%
\end{equation*}
\end{enumerate}
\end{remarks}

\begin{notation}
For $N\in \mathbb{S}_{B}$ and $L\in \mathbb{S}^{B},$ we have the following
morphisms in $\mathbb{S}_{B}^{B}:$%
\begin{eqnarray*}
\gamma _{N} &:&N\otimes _{S}B\longrightarrow N\otimes _{S}^{a}B,\text{ }%
n\otimes _{S}b\mapsto (n\otimes _{S}1_{B})\Delta (b)=\sum nb_{1}\otimes
_{S}b_{2}; \\
\gamma ^{L} &:&L\otimes _{S}^{c}B\longrightarrow L\otimes _{S}B,\text{ }%
l\otimes _{S}^{c}b\mapsto \rho ^{L}(l)(1_{B}\otimes _{S}b)=\sum
l_{<0>}\otimes _{S}l_{<1>}b.
\end{eqnarray*}%
In particular, $\gamma _{B}:B\otimes _{S}B\longrightarrow B\otimes _{S}^{a}B$
and $\gamma ^{B}:B\otimes _{S}^{c}B\longrightarrow B\otimes _{S}B$ are
morphisms in $\mathbb{S}_{B}^{B}.$
\end{notation}

\begin{definition}
Let $\mathfrak{A}$ be a category with finite limits and finite colimits. A
functor $F:\mathfrak{A}\longrightarrow \mathfrak{B}$ is said to be \emph{%
left-exact}\textbf{\ }(\emph{right-exact}\textbf{)} iff $F$ preserves finite
limits (finite colimits). Moreover, $F$ is called \emph{exact} iff $F$ is
left-exact and right-exact.
\end{definition}

The following technical result will be needed in the sequel; the proof is
straightforward.

\begin{lemma}
\label{NotB}Let $(B,\mu ,\eta ,\Delta ,\varepsilon )$ be an $S$%
-bisemialgebra.

\begin{enumerate}
\item $(\mathcal{F},-\otimes _{S}^{a}B)$ is an adjoint pair of functors,
where $\mathcal{F}:\mathbb{S}_{B}^{B}\longrightarrow \mathbb{S}_{B}$ is the
forgetful functor. Consequently, $-\otimes _{S}^{a}B:\mathbb{S}%
_{B}\longrightarrow \mathbb{S}_{B}^{B}$ is left exact and preserves all
limits \emph{(e.g.} equalizers, kernels and direct products\emph{)}.

\item $(-\otimes _{S}^{c}B,\mathcal{G)}$ is an adjoint pair of functors,
where $\mathcal{G}:\mathbb{S}_{B}^{B}\longrightarrow \mathbb{S}^{B}$ is the
forgetful functor. Consequently, $-\otimes _{S}^{c}B:\mathbb{S}%
^{B}\longrightarrow \mathbb{S}_{B}^{B}$ is right exact and preserves all
colimits \emph{(e.g.} coequalizers, cokernels and direct coproducts\emph{)}.
\end{enumerate}
\end{lemma}

\begin{Beweis}
It is a well-known fact that a left (right) adjoint functor is right (left)
exact and preserves all colimits (limits).

\begin{enumerate}
\item For every $M\in \mathbb{S}_{B}^{B}$ and $N_{B},$ we have a natural
isomorphism of $S$-semimodules%
\begin{equation}
\mathrm{Hom}_{B}^{B}(M,N\otimes _{S}^{a}B)\longrightarrow \mathrm{Hom}_{B}(%
\mathcal{F}(M),N),\text{ }f\mapsto (\vartheta _{N}^{r}\circ (N\otimes
_{S}\varepsilon ))\circ f  \label{N_B}
\end{equation}%
with inverse $h\mapsto \lbrack h\otimes _{S}B\circ \rho ^{M}].$

\item For every $M\in \mathbb{S}_{B}^{B}$ and $N^{B},$ we have a natural
isomorphism of $S$-semimodules%
\begin{equation}
\mathrm{Hom}_{B}^{B}(N\otimes _{S}^{c}B,M)\longrightarrow \mathrm{Hom}^{B}(N,%
\mathcal{G}(M)),\text{ }g\mapsto g(-\otimes _{S}1_{B})  \label{N^B}
\end{equation}%
with inverse $h\mapsto \lbrack \rho _{M}\circ (h\otimes _{S}B)].$
\end{enumerate}
\end{Beweis}

\subsection*{Integrals}

\qquad As before, we let $(B,\mu ,\eta ,\Delta ,\varepsilon )$ be an $S$%
-bisemialgebra.

\begin{punto}
A \emph{left} (\emph{total}) \emph{integral \textbf{on} }$B$ is a left $B$%
-colinear morphism $t\in B^{\ast }:=\mathrm{Hom}_{S}(B,S)$ (with $%
t(1_{B})=1_{S}$), equivalently $\sum b_{1}t(b_{2})=1_{B}t(b)$ for every $%
b\in B$ (and $t(1_{B})=1_{S}$). The \emph{right }(\emph{total}) \emph{%
integrals} are defined symmetrically. For a justification of the
terminology, we refer to \cite{Mon1993} (see also \cite[p. 181]{DNR2001}).
With $\int^{l}\leq B^{\ast }$ ($\int^{l,1}\leq B^{\ast }$) we denote the $S$%
-subsemimodule of (total) left integrals on $B;$ symmetrically, we denote
with $\int^{r}\leq B^{\ast }$ ($\int^{r,1}\leq B^{\ast }$) the $S$%
-subsemimodule of (total) right integrals on $B.$
\end{punto}

\begin{definition}
Consider $B^{\ast }$ as a $(B,B)$-bisemimodule in the canonical way. If $%
_{S}B$ is an $\alpha $-semimodule, then we call $\mathrm{Rat}^{B}(_{B^{\ast
}}B^{\ast })$ ($^{B}\mathrm{Rat}(B_{B^{\ast }}^{\ast })$) the \emph{left} (%
\emph{right}) \emph{trace ideal} of $B^{\ast }.$
\end{definition}

\begin{lemma}
\label{left-int}Let $t\in B^{\ast }:=\mathrm{Hom}_{S}(B,S).$

\begin{enumerate}
\item Assume that $B$ is $S$-cogenerated in $\mathbb{S}_{S}.$

\begin{enumerate}
\item $t\in \int^{l}$ if and only if $f\ast t=f(1_{B})t$ for every $f\in
B^{\ast }.$

\item $\int^{l}$ is an ideal of $B^{\ast }.$
\end{enumerate}

\item Let $_{S}B$ be an $\alpha $-semimodule and $B^{\ast rat}:=\mathrm{Rat}%
^{B}(_{B^{\ast }}B^{\ast }).$ Then $t\in \int^{l}$ if and only if $\rho
^{B^{\ast rat}}(t)=t\otimes _{S}1_{B}.$
\end{enumerate}
\end{lemma}

\begin{ex}
Consider the polynomial $S$-bisemialgebra $S[x]$ with $\Delta (x)=x\otimes
_{S}x$ and $\varepsilon (x)=1_{S}.$ It is clear that%
\begin{equation*}
t:S[x]\longrightarrow S,\text{ }p(x)\mapsto \delta _{1,p(x)}
\end{equation*}%
is a \emph{total} left (right) integral on $S[x]:$ let $f\in S[x]^{\ast }.$
For every $n\in \mathbb{N}_{0},$ we have 
\begin{equation*}
(f\ast t)(x^{n})=f(x^{n})t(x^{n})=\delta _{1,x^{n}}f(x^{n})=\delta
_{0,n}f(x^{n})=f(1)\delta _{0,n}=f(1)\delta _{1,x^{n}}.
\end{equation*}%
Since $S[x]\simeq S^{(\mathbb{N}_{0})}\hookrightarrow S^{\mathbb{N}_{0}},$
it is $S$-cogenerated as an $S$-semimodule, we conclude that $t$ is a left
integral on $(S[x],\Delta ,\varepsilon ).$ Similarly, $t$ is a right
integral on $(S[x],\Delta ,\varepsilon ).$
\end{ex}

\begin{ex}
Consider the bisemialgebra $B=\mathbb{B}[x],$ where $B$ is the \emph{Boolean
semifield}, with%
\begin{eqnarray*}
\Delta (x) &=&x\otimes _{\mathbb{B}}1+1\otimes _{\mathbb{B}}x,\text{ }\Delta
(1)=1\otimes _{S}1; \\
\varepsilon (x) &=&0,\text{ }\varepsilon (1)=1.
\end{eqnarray*}%
Let $t\in B^{\ast }$ be a left integral on $B.$ We have $%
1_{B}t(x)=xt(1)+1_{B}t(x),$ whence $t(1)=0.$ We prove \emph{by induction}
that $t(x^{n})=0$ for $n\geq 1.$ Since $\Delta (x^{2})=x^{2}\otimes _{%
\mathbb{B}}1_{B}+x\otimes _{\mathbb{B}}x+1_{B}\otimes _{\mathbb{B}}x^{2},$
it follows that $1_{B}t(x^{2})=xt(x)+1_{B}t(x),$ whence $t(x)=0.$ Consider $%
n\geq 2$ and assume that $t(x^{n-i})=0$ for all $i=1,\cdots ,n.$ We have $%
\Delta (x^{n+1})=x^{n+1}\otimes _{\mathbb{B}}1+x^{n}\otimes _{\mathbb{B}%
}x+\cdots +x\otimes _{\mathbb{B}}x^{n}+1\otimes _{\mathbb{B}}x^{n+1}$ and so 
$1_{B}t(x^{n+1})=xt(x^{n})+1_{B}t(x^{n+1}),$ whence $t(x^{n})=0.$
Consequently, $\int^{l}B=0.$ Similarly, we can prove that $\int^{r}B=0.$
\end{ex}

\begin{punto}
Let $M$ be a left $B$-semimodule. The set of \emph{invariants} of $M$ is%
\begin{equation*}
^{B}M:=\{m\in M\mid bm=\varepsilon (b)m\text{ for every }b\in B\}.
\end{equation*}%
Moreover, we have an isomorphism of $S$-semimodules%
\begin{equation*}
\mathrm{Hom}_{B-}(S,M)\longrightarrow \text{ }^{B}M,\text{ }f\longmapsto
f(1_{S}).
\end{equation*}
\end{punto}

\begin{definition}
A \emph{left integral \textbf{in} }$B$ is an invariant of $_{B}B,$ \emph{i.e.%
} an element of%
\begin{equation*}
^{B}B:=\{\varpi \in B\mid b\varpi =\varepsilon (b)\varpi \text{ for every }%
\varpi \in B\}.
\end{equation*}%
We say that a left integral $\varpi $ in $H$ is \emph{normalized} iff $%
\varepsilon (\varpi )=1_{S}.$ Symmetrically, one can define (\emph{normalized%
}) \emph{right integrals in }$B.$ With $\int_{l}B$ ($\int_{l,1}B$) we denote
the set of (normalized) left integral in B and with $\int_{r}B$ ($%
\int_{r,1}B $) we denote the set of (normalized) right integrals in $B.$
\end{definition}

\begin{ex}
Let $G=\{g_{1},\cdots ,g_{n}\}$ be a finite group and consider the $S$%
-bisemialgebra $B=S[G]$ with $\Delta (g_{i})=g_{i}\otimes _{S}g_{i}$ and $%
\varepsilon (g_{i})=1_{S}$ for $i=1,\cdots ,n.$ Then $\varpi :=g_{1}+\cdots
+g_{n}$ is a left (right) integral in $B:$ indeed, for every $g\in G,$ we
have 
\begin{equation*}
g\varpi
=g(\sum_{i=1}^{n}g_{i})=\sum_{i=1}^{n}(gg_{i})=\sum_{j=1}^{n}g_{j}=%
\varepsilon (g)\varpi .
\end{equation*}
\end{ex}

\begin{punto}
Let $M$ be a right $B$-semicomodule. The set of \emph{coinvariants} of $M$ is%
\begin{equation*}
M^{\mathrm{co}B}:=\mathrm{Eq}(\rho ^{M},g)=\{m\in M\mid \rho
^{M}(m)=m\otimes _{S}1_{B}\},\text{ where }g(m):=m\otimes _{S}1_{B}
\end{equation*}%
\begin{equation*}
\xymatrix{ M^{\mathrm{co}B} \ar[r] & M \ar@<1ex>^{\rho ^M}[rr] \ar_{g}[rr] &
& M \otimes_S B}
\end{equation*}
\end{punto}

\begin{lemma}
\label{coB-prop}

\begin{enumerate}
\item For all $L\in \mathbb{S}_{S}$ and $M\in \mathbb{S}^{B},$ we have an
isomorphism%
\begin{equation}
\mathrm{Hom}^{B}(L,M)\simeq \mathrm{Hom}_{S}(L,M^{\mathrm{co}B}).
\end{equation}

\item We have an isomorphism of $S$-semimodules%
\begin{equation}
\mathrm{Hom}^{B}(S,M)\longrightarrow M^{\mathrm{co}B},\text{ }f\longmapsto
f(1_{S}).
\end{equation}

\item $B^{\mathrm{co}B}=S1_{B}.$
\end{enumerate}
\end{lemma}

\begin{proposition}
\label{co-inv}Let $_{S}B$ be an $\alpha $-semimodule and consider the trace
ideal $B^{\ast rat}:=\mathrm{Rat}^{B}(_{B^{\ast }}B^{\ast }).$

\begin{enumerate}
\item For every $M\in \mathbb{S}^{B},$ we have $M^{\mathrm{co}B}=$ $%
^{B^{\ast }}M.$

\item $(B^{\ast rat})\mathrm{^{\mathrm{co}B}}=$ $^{B^{\ast }}(B^{\ast rat}).$

\item If $_{S}B$ is finitely generated and projective, then $(B^{\ast })^{%
\mathrm{co}B}=$ $^{B^{\ast }}B^{\ast }.$
\end{enumerate}
\end{proposition}

\begin{Beweis}
$(1)$ The proof is similar to that of \cite[14.13 (1)]{BW2003}.

$(2)$ Set $M=B^{\ast rat}$ in (1) and notice that $B^{\ast rat}\in \mathbb{S}%
^{B}$ by (\ref{rat-co}).

$(3)$ Since $_{S}B$ is finitely generated and projective, $B\simeq B^{\ast
\ast }$ and $B^{\ast }$ is finitely generated and projective. It follows
that we have a canonical isomorphism of $S$-semimodules $B^{\ast }\otimes
_{S}B^{\ast \ast }\overset{\beta _{B^{\ast }}}{\longrightarrow }\mathrm{Hom}%
_{S}(B^{\ast },B^{\ast })$ \cite[Proposition 3.7]{KN2011}. Notice that $%
\alpha _{B^{\ast }}^{B}$ is in fact the composition of the following
isomorphisms $B^{\ast }\otimes _{S}B\simeq B^{\ast }\otimes _{S}B^{\ast \ast
}\simeq \mathrm{Hom}_{S}(B^{\ast },B^{\ast }),$ whence%
\begin{equation*}
B^{\ast rat}:=\mathrm{Rat}^{B}(_{B^{\ast }}B^{\ast })=\widetilde{\rho }%
_{B^{\ast }}^{-1}(\alpha _{B^{\ast }}^{B}(B^{\ast }\otimes _{S}B))=%
\widetilde{\rho }_{B^{\ast }}^{-1}(\mathrm{Hom}_{S}(B^{\ast },B^{\ast
}))=B^{\ast }
\end{equation*}%
and the result follows from (2).$\blacksquare $
\end{Beweis}

\begin{proposition}
\label{coB}

\begin{enumerate}
\item $(-\otimes _{S}B,(-)^{\mathrm{co}B})$ is an adjoint pair of functors.

\item $(-\otimes _{S}B,\mathrm{Hom}_{B}^{B}(B,-))$ is an adjoint pair of
functors.

\item We have a natural isomorphism of functors 
\begin{equation}
\mathrm{Hom}_{B}^{B}(B,-)\simeq (-)^{\mathrm{co}B}.  \label{BB=coB}
\end{equation}

\item $\mathrm{Hom}_{B}^{B}(B,-)$ is left exact and preserves all limits 
\emph{(e.g.} equalizers, kernels and direct products\emph{)}.

\item We have an isomorphism of semirings%
\begin{equation}
\mathrm{End}_{B}^{B}(B)\simeq B^{\mathrm{co}B}=B1_{S}.
\end{equation}
\end{enumerate}
\end{proposition}

\begin{Beweis}
It is clear that $(-)^{\mathrm{co}B}:\mathbb{S}_{B}^{B}\longrightarrow 
\mathbb{S}_{S}$ and $\mathrm{Hom}_{B}^{B}(B,-):\mathbb{S}_{B}^{B}%
\longrightarrow \mathbb{S}_{S}$ are functors.

\begin{enumerate}
\item For every $M\in \mathbb{S}_{S}$ and $N\in \mathbb{S}_{B}^{B},$ we have
a natural isomorphism%
\begin{equation}
\mathrm{Hom}_{B}^{B}(M\otimes _{S}B,N)\simeq \mathrm{Hom}_{S}(M,N^{\mathrm{co%
}B}),\text{ }f\mapsto \lbrack m\mapsto f(m\otimes _{S}1_{B})]
\label{ot-B-co}
\end{equation}%
with inverse $g\mapsto \lbrack m\otimes _{S}b\mapsto g(m)b].$

\item For every $M\in \mathbb{S}_{S}$ and $N\in \mathbb{S}_{B}^{B},$ we have
a natural isomorphism%
\begin{equation}
\mathrm{Hom}_{B}^{B}(M\otimes _{S}B,N)\simeq \mathrm{Hom}_{S}(M,\mathrm{Hom}%
_{B}^{B}(B,N)),\text{ }f\mapsto \lbrack -\mapsto f(-\otimes _{S}1_{B})]
\label{ot-B-BB}
\end{equation}%
with inverse $g\mapsto \lbrack m\otimes _{S}b\mapsto g(m)(b)].$

\item This follows from (1) and (2) by the uniqueness of the right adjoint
of the functor $-\otimes _{S}B:\mathbb{S}_{S}\longrightarrow \mathbb{S}%
_{B}^{B}.$ In fact, substituting $M=S$ in (\ref{ot-B-co}) yields a natural
isomorphism for every $N\in \mathbb{S}_{B}^{B}:$%
\begin{equation}
\mathrm{Hom}_{B}^{B}(B,N)\simeq N^{\mathrm{co}B},\text{ }f\mapsto f(1_{B})
\label{f(1)}
\end{equation}%
with inverse $n\mapsto \lbrack b\mapsto nb].$

\item This follows directly from the fact that $\mathrm{Hom}_{B}^{B}(B,-)$
has a left adjoint by (2).

\item Set $N:=B$ in (2). It is clear that the isomorphism obtained is in
fact a morphism of semirings.$\blacksquare $
\end{enumerate}
\end{Beweis}

We present now the main reconstruction result in this section:

\begin{theorem}
\label{bimonad}

\begin{enumerate}
\item We have an isomorphism of categories $\mathbf{SBialg}_{S}\simeq 
\mathbf{Bimonoid}(\mathbb{S}_{S}).$

\item Let $B$ be an $S$-semimodule. There is a bijective correspondence
between the structures of $S$-bisemialgebras on $B,$ the bimonad structures
on $B\otimes _{S}-:\mathbb{S}_{S}\longrightarrow \mathbb{S}_{S}$ and the
bimonad structures on $-\otimes _{S}B:\mathbb{S}_{S}\longrightarrow \mathbb{S%
}_{S}.$

\item We have isomorphisms of categories%
\begin{eqnarray*}
\mathbb{S}_{B}^{B} &\simeq &(\mathbb{S}_{B})^{-\otimes _{S}B}\simeq ((%
\mathbb{S}_{S})_{-\otimes _{S}B})^{-\otimes _{S}B}; \\
&\simeq &(\mathbb{S}^{B})_{-\otimes _{S}B}\simeq ((\mathbb{S}_{S})^{-\otimes
_{S}B})_{-\otimes _{S}B}.
\end{eqnarray*}
\end{enumerate}
\end{theorem}

\begin{Beweis}
(1) and (3) follow directly from the definitions. The proof of (2) is along
the liens of that of \cite[Theorem 3.9]{Ver} taking into consideration that $%
\mathbb{S}_{S}$ is cocomplete, that $S$ is a regular generator in $\mathbb{S}%
_{S}$ \cite{Gol1999} and the fact that $-\otimes _{S}X\simeq X\otimes _{S}-:%
\mathbb{S}_{S}\longrightarrow \mathbb{S}_{S}$ preserves colimits for every $%
S $-semimodule $X.\blacksquare $
\end{Beweis}

\subsection*{An Application}

\begin{punto}
(\cite{Wor2012}) A \emph{right }$S$\emph{-linear automaton} is a datum $%
\mathbf{A}=(M,A,\mathbf{s},\mathbf{\rho },\Omega ),$ where $A$ is an $S$%
-semialgebra, $M$ is a right $A$-semimodule, $\mathbf{s}\in M$ (called a 
\emph{starting vector}) and $\Omega \in \mathrm{Hom}_{S}(M,S)$ (called an 
\emph{observation function}). The \emph{language accepted by a right }$S$%
\emph{-linear automaton} $A$ is the $S$-linear map $\mathbf{\rho }%
:A\longrightarrow S,$ $a\mapsto \Omega (\mathbf{s}a).$
\end{punto}

\begin{punto}
Let $B$ be an $S$-bisemialgebra. If $\mathbf{A}=(M,B,\mathbf{s},\rho ,\Omega
)$ and $\mathbf{A}^{\prime }=(M^{\prime },B,\mathbf{s}^{\prime },\rho
^{\prime },\Omega ^{\prime })$ are two left $S$-linear automata, then $%
\mathbf{A}\otimes _{S}\mathbf{A}^{\prime }:=(M\otimes _{S}^{b}N,B,\mathbf{s}%
\otimes _{S}\mathbf{s}^{\prime },\rho _{M\otimes _{S}^{b}N},\Omega \otimes
_{S}\Omega ^{\prime })$ is a right $S$-linear automaton and the language
accepted by $A\otimes _{S}A^{\prime }$ is $\mathbf{\rho }_{A\otimes
_{S}A^{\prime }}:=\mathbf{\rho }\ast \mathbf{\rho }^{\prime }.$
\end{punto}

\section{Doi-Koppinen Semimodules}

The class of Doi-Koppinen modules over rings was introduced independently by
Y. Doi \cite{Doi1992} and M. Koppinen \cite{Kop1995}. In this section, we
extend these notions and some results on them to Doi-Koppinen semimodules
over semirings. Throughout this section, $(B,\mu ,\eta ,\Delta ,\varepsilon
) $ is an $S$-bisemialgebra.

\begin{definition}
\begin{enumerate}
\item A \emph{right }$B$\emph{-semimodule semialgebra} is an $S$-semialgebra 
$(A,\mu _{A},\eta _{A})$ with a right $B$-semimodule structure such that $%
\mu _{A}$ and $\eta _{A}$ are $B$\emph{-linear}, i.e. 
\begin{equation}
(a\widetilde{a})b=\sum (ab_{1})(\widetilde{a}b_{2})\text{ and }%
1_{A}b=\varepsilon _{B}(b)1_{A}\text{ for all }a,\widetilde{a}\in A\text{
and }b\in B.  \label{rma}
\end{equation}%
Symmetrically, one defines a \emph{left }$B$\emph{-semimodule algebra}.

\item A \emph{right }$B$\emph{-semimodule semicoalgebra }is an $S$%
-semicoalgebra $(C,\Delta _{C},\varepsilon _{C})$ with a right $B$%
-semimodule structure such that $\Delta _{C}$ and $\varepsilon _{C}$ are $B$%
\emph{-linear}, i.e. 
\begin{equation}
\sum (cb)_{1}\otimes _{S}(cb)_{2}=\sum c_{1}b_{1}\otimes c_{2}b_{2}\text{
and }\varepsilon _{C}(cb)=\varepsilon _{C}(c)\varepsilon _{B}(b)\text{ for
all }c\in C\text{ and }b\in B.  \label{rmc}
\end{equation}%
Symmetrically, one defines a \emph{left }$B$\emph{-semimodule semicoalgebra}.

\item A \emph{right }$B$\emph{-semicomodule semialgebra} is an $S$%
-semialgebra $(A,\mu _{A},\eta _{A})$ with a right $B$-semicomodule
structure such that $\mu _{A}$ and $\eta _{A}$ are $B$\emph{-colinear}, i.e.%
\begin{equation}
\sum (ab)_{<0>}\otimes _{S}(ab)_{<1>}=\sum a_{<0>}b_{<0>}\otimes
_{S}a_{<1>}b_{<1>}\text{ and }\sum 1_{<0>}\otimes _{S}1_{<1>}=1_{A}\otimes
1_{B}.  \label{rca}
\end{equation}%
Symmetrically, one defines a \emph{left }$B$\emph{-semicomodule semialgebra}.

\item A \emph{right }$B$\emph{-semicomodule semicoalgebra} is an $S$%
-semicoalgebra $(C,\Delta _{C},\varepsilon _{C})$ with a right $B$%
-semicomodule structure such that $\Delta _{C}$ and $\varepsilon _{C}$ are $%
B $\emph{-colinear}, i.e. 
\begin{equation}
\sum c_{<0>1}\otimes c_{<0>2}\otimes c_{<1>}=\sum c_{1<0>}\otimes
c_{2<0>}\otimes c_{1<1>}c_{2<1>}\text{ and }\sum \varepsilon
_{C}(c_{<0>})c_{<1>}=\varepsilon _{C}(c)1_{B}.  \label{com-coal}
\end{equation}%
Symmetrically, one defines a \emph{left }$B$\emph{-semicomodule semicoalgebra%
}.
\end{enumerate}
\end{definition}

\begin{punto}
\label{DK}A \emph{right-right Doi-Koppinen structure} over $S$ is a triple $%
(B,A,C)$ consisting of an $S$-bisemialgebra $B,$ a right $B$-semicomodule
semialgebra $A$ and a right $B$-semimodule semicoalgebra $C.$ A \emph{%
right-right Doi-Koppinen semimodule} for $(B,A,C)$ is a right $A$-semimodule 
$M,$ which is also a right $C$-semicomodule such that 
\begin{equation*}
\sum (ma)_{<0>}\otimes _{S}(ma)_{<1>}=\sum m_{<0>}a_{<0>}\otimes
_{S}m_{<1>}a_{<1>}\text{ for all }m\in M\text{ and }a\in A.
\end{equation*}%
With $\mathbb{S}(B)_{A}^{C}$ we denote the category of right-right
Doi-Koppinen semimodules and $A$-linear $C$-colinear morphisms.
\end{punto}

The following result is easy to prove.

\begin{lemma}
Let $(B,A,C)$ be a right-right Doi-Koppinen structure over $S.$

\begin{enumerate}
\item $\#^{op}(C,A):=\mathrm{Hom}_{S}(C,A)$ is an $A$-semiring with $(A,A)$%
-bisemimodule structure%
\begin{equation*}
(af)(c):=\sum a_{<0>}f(ca_{<1>})\text{ and }(fa)(c):=f(c)a,
\end{equation*}%
multiplication%
\begin{equation}
(f\cdot g)(c)=\sum f(c_{2})_{<0>}g(c_{1}f\left( c_{2}\right) _{<1>})
\label{DK-mult}
\end{equation}%
and unity $\eta _{A}\circ \varepsilon _{C}.$

\item $\mathcal{C}:=A\otimes _{S}C$ is an $A$-semicoring and $%
\#^{op}(C,A)\simeq $ $^{\ast }\mathcal{C}$ as $A$-semirings.
\end{enumerate}
\end{lemma}

\begin{exs}
\begin{enumerate}
\item $C=B$ is a right $B$-semimodule semicoalgebra with structure map $\mu
_{B}$ and so $(B,A,B)$ is a right-right Doi-Koppinen structure for every
right $B$-semicomodule semialgebra $A.$ In this case, $\mathbb{S}(B)_{A}^{B}=%
\mathbb{S}_{A}^{B},$ the category of \emph{relative Hopf semimodules} (cf. 
\cite{Doi1983}).

\item $A=B$ is a right $B$-semicomodule semialgebra with structure map $%
\Delta _{B}$ and so $(B,B,C)$ is a right-right Doi-Koppinen structure for
every right $B$-semimodule semicoalgebra $C.$ In this case, $\mathbb{S}%
(B)_{B}^{C}=\mathbb{S}_{[C,B]},$ the category of \emph{Doi's }$[C,B]$\emph{%
-semimodules }(cf.\emph{\ }\cite{Doi1983}).

\item Setting $A=B=C,$ we notice that $(B,B,B)$ is a right-right
Doi-Koppinen structure and that $\mathbb{S}(B)_{B}^{B}=\mathbb{S}_{B}^{B},$
the category of \emph{Hopf semimodules} (cf. \cite[4.1]{Swe1969}).
\end{enumerate}
\end{exs}

\qquad The following result is easy to prove.

\begin{lemma}
\label{rr-ring}Let $(B,A,C)\;$be a right-right Doi-Koppinen structure over $%
S $ and consider the corresponding $A$-semicoring $\mathcal{C}:=A\otimes
_{S}C. $

\begin{enumerate}
\item $A\#^{op}C^{\ast }:=A\otimes _{S}C^{\ast }$ is an $A$-semiring with $%
(A,A)$-bisemimodule structure 
\begin{equation}
\widetilde{a}(a\#f):=\sum \widetilde{a}_{<0>}a\#\widetilde{a}_{<1>}f\text{
and }(a\#f)\widetilde{a}:=a\widetilde{a}\#f,  \label{AC*-bimod}
\end{equation}%
multiplication%
\begin{equation}
(a\#f)\cdot (b\#g):=\sum a_{<0>}b\#(a_{<1>}g)\ast f  \label{op-smash}
\end{equation}%
and unity $1_{A}\#\varepsilon _{C}.$ Moreover,%
\begin{equation*}
\eta :A\longrightarrow A\#^{op}C^{\ast },a\mapsto a\#\varepsilon _{C}
\end{equation*}%
is a morphism of $A$-semirings.

\item $P:=(A\#^{op}C^{\ast },\mathcal{C})$ is a measuring left $A$-pairing 
\emph{(}in the sense of \cite{Abu-c}\emph{)}.
\end{enumerate}
\end{lemma}

\begin{Beweis}
\begin{enumerate}
\item Clear.

\item It is clear that%
\begin{equation*}
\kappa _{P}:A\#^{op}C^{\ast }\longrightarrow \text{ }^{\ast }\mathcal{C},%
\text{ }a\#f\mapsto \lbrack \widetilde{a}\otimes _{S}c\mapsto \widetilde{a}%
af(c)]
\end{equation*}%
is a morphism of $A$-semirings.$\blacksquare $
\end{enumerate}
\end{Beweis}

\begin{theorem}
\label{HAC-iso}Let $(B,A,C)$ be a right-right Doi-Koppinen structure and
consider the corresponding $A$-semicoring $\mathcal{C}:=A\otimes _{S}C.$ We
have an isomorphism of categories%
\begin{equation*}
\mathbb{S}(B)_{A}^{C}\simeq \mathbb{S}^{\mathcal{C}}.
\end{equation*}
\end{theorem}

\begin{Beweis}
It can be shown that $(A,C,\psi )$ is a \emph{right-right entwining structure%
} in the symmetric monoidal category $(\mathbb{S}_{S},\otimes _{S},S;\mathbf{%
\tau }),$ where%
\begin{equation}
\psi :C\otimes _{A}A\longrightarrow A\otimes _{A}C,\text{ }c\otimes
_{A}a\mapsto \sum a_{<0>}\otimes _{A}ca_{<1>}.  \label{psi-DK}
\end{equation}%
By arguments similar to those in \cite{Brz1999} (see also \cite[Ch. 5]%
{BW2003}), one can show that $\mathbb{S}(B)_{A}^{C}\simeq \mathbb{S}%
_{A}^{C}(\psi )\simeq \mathbb{S}^{\mathcal{C}},$ where $\mathbb{S}%
_{A}^{C}(\psi )$ is the associated category of \emph{right-right entwined
semimodules}.$\blacksquare $
\end{Beweis}

\begin{proposition}
\label{Porst-AC}Let $(B,A,C)$ be a right-right Doi-Koppinen structure and
the associated category $\mathbb{S}(B)_{A}^{C}$ of right-right Doi-Koppinen
semimodules.

\begin{enumerate}
\item $\mathbb{S}(B)_{A}^{C}$ is comonadic, locally presentable and a
covariety \emph{(}in the sense of \cite{AP2003}\emph{)}.

\item The forgetful functor $\mathcal{F}:\mathbb{S}(B)_{A}^{C}%
\longrightarrow \mathbb{S}_{A}$ creates all colimits and isomorphisms.

\item $\mathbb{S}(B)_{A}^{C}$ is cocomplete, \emph{i.e.} $\mathbb{S}%
(B)_{A}^{C}$ has all \emph{(}small\emph{)} colimits, \emph{e.g.}
coequalizers, cokernels, pushouts, directed colimits and direct sums.
Moreover, the colimits are formed in $\mathbb{S}_{A}.$

\item $\mathbb{S}(B)_{A}^{C}$ is complete, \emph{i.e.} $\mathbb{S}%
(B)_{A}^{C} $ has all \emph{(}small\emph{)} limits, \emph{e.g.} equalizers,
kernels, pullbacks, inverse limits and direct products. Moreover, the
forgetful functor $\mathcal{F}$ creates all limits preserved by $-\otimes
_{A}(A\otimes _{S}C):\mathbb{S}_{A}\longrightarrow \mathbb{S}_{A}.$
\end{enumerate}
\end{proposition}

\begin{Beweis}
The result is an immediate consequence of \cite[Proposition 2.22]{Abu-c}
taken into consideration that $\mathcal{C}:=A\otimes _{S}C$ is an $A$%
-semicoring and that $\mathbb{S}(B)_{A}^{C}\simeq \mathbb{S}^{\mathcal{C}%
}.\blacksquare $
\end{Beweis}

\begin{remark}
Let $(B,A,C)\;$be a right-right Doi-Koppinen structure over $S.$ While it is
guaranteed that the category $\mathbb{S}(B)_{A}^{C}$ has kernels, these are
not necessarily formed in $\mathbb{S}_{A}.$ Indeed, if $_{S}C$ is flat, then 
$_{A}\mathcal{C}$ is flat and it follows by \cite[Proposition 2.26]{Abu-c}
that all equalizers in $\mathbb{S}(B)_{A}^{C}\simeq \mathbb{S}^{\mathcal{C}}$
are formed in $\mathbb{S}_{A}.$
\end{remark}

\begin{theorem}
\label{AC-sigma}Let $(B,A,C)$ be a right-right Doi-Koppinen structure and
consider the corresponding $A$-semicoring $\mathcal{C}:=A\otimes _{S}C.$ If $%
P:=(A\#^{op}C^{\ast },\mathcal{C})$ satisfies the $\alpha $-condition, then
we have an isomorphisms of categories%
\begin{equation*}
\mathbb{S}(B)_{A}^{C}\simeq \mathrm{Rat}^{\mathcal{C}}(\mathbb{S}%
_{A\#^{op}C^{\ast }})=\sigma \lbrack \mathcal{C}_{A\#^{op}C^{\ast }}].
\end{equation*}
\end{theorem}

\begin{Beweis}
The isomorphism $\mathbb{S}^{\mathcal{C}}\simeq \mathrm{Rat}^{\mathcal{C}}(%
\mathbb{S}_{A\#^{op}C^{\ast }})$ follows by \cite[Theorem 3.16]{Abu-c}. A
similar argument to that of \cite[Theorem 3.22]{Abu-c} shows that $\mathbb{S}%
^{\mathcal{C}}\simeq \mathrm{Rat}^{\mathcal{C}}(\mathbb{S}_{B\#^{op}B^{\ast
}})=\sigma \lbrack \mathcal{C}_{B\#^{op}B^{\ast }}].\blacksquare $
\end{Beweis}

We provide now an example in which the assumption of Theorem \ref{HAC-iso}
(2) holds:

\begin{ex}
\label{ex-C-C}Let $(B,A,C)$ be a right-right Doi-Koppinen structure with $%
_{S}C$ an $\alpha $-semimodule. Consider the left measuring $A$-pairing $%
P:=(A\#^{op}C^{\ast },\mathcal{C})$ and let $\phi :C^{\ast }\longrightarrow
A\otimes _{S}C^{\ast },$ $f\mapsto 1_{A}\otimes _{S}f.$ For every right $A$%
-semimodule $M,$ we have the following commutative diagram%
\begin{equation*}
\xymatrix{ M \otimes_A (A \otimes_S C) \ar[rr]^{\alpha_M ^{P}} \ar@{=}[d] &
& {\rm Hom}_{-A} (A \otimes_S C^*,M) \ar[d]^{(\phi,M)} \\ M \otimes_S C
\ar[rr]_{\alpha_M ^{C}} & & {\rm Hom}_{S} (C^*,M) }
\end{equation*}%
Since $\alpha _{M}^{C}$ is injective, it follows that $\alpha _{M}^{P}$ is
injective. We claim that $\alpha _{M}^{P}$ is uniform for every $M_{A}.$ Let 
$M$ be a right $A$-semimodule and let $h\in \overline{\alpha
_{M}^{P}(M\otimes _{A}(A\otimes _{S}C)},$ \emph{i.e. }$h+\alpha ^{P}(\sum
m_{i}\otimes _{A}(a_{i}\otimes _{S}c_{i}))=\alpha ^{P}(\sum m_{i}^{\prime
}\otimes _{A}(a_{i}^{\prime }\otimes c_{i}^{\prime })).$ For every $g\in
C^{\ast }$ we have%
\begin{eqnarray*}
h(1_{A}\otimes _{S}g)+\alpha _{M}^{P}(\sum m\otimes _{A}(a_{i}\otimes
_{S}c_{i}))(1_{A}\otimes _{S}g) &=&\alpha _{M}^{P}(\sum m_{i}\otimes
_{A}(a_{i}\otimes _{S}c_{i}))(1_{A}\otimes _{S}g) \\
(\phi ,M)(h)(g)+\alpha _{M}^{C}(\sum m_{i}a_{i}\otimes _{S}c_{i})(g)
&=&\alpha _{M}^{C}(\sum m_{i}a_{i}\otimes _{S}c_{i})(g),
\end{eqnarray*}%
whence $(\phi ,M)(h)\in \overline{\alpha _{M}^{C}(M\otimes _{S}C)}.$ Since $%
\alpha _{M}^{C}$ is uniform (by our assumption on $_{S}C$), there exists $%
\sum m_{j}^{\prime }\otimes _{S}c_{j}^{\prime }\in M\otimes _{S}C$ $(\varphi
,M)(h)$ such that for every $g\in C^{\ast }$ we have 
\begin{equation*}
h(1\otimes g)=(\varphi ,M)(h)(g)=\alpha ^{C}(\sum m_{j}\otimes
c_{j})(g)=\sum m_{j}g(c_{j}).
\end{equation*}%
Noting that $h$ is right $A$-linear, it follows that for every $%
\sum_{l}a_{l}^{\prime }\otimes g_{l}^{\prime }\in A\otimes _{S}C^{\ast }$ we
have%
\begin{equation*}
\begin{tabular}{lllll}
$h(\sum_{l}a_{l}^{\prime }\otimes g_{l}^{\prime })$ & $=$ & $%
h(\sum_{l}(1_{A}\otimes _{S}g_{l}^{\prime })a_{l}^{\prime })$ & $=$ & $%
\sum_{l}h(1_{A}\otimes _{S}g_{l}^{\prime })a_{l}^{\prime }$ \\ 
& $=$ & $\sum_{k}\sum_{j}m_{j}a_{k}^{\prime }g_{k}^{\prime }(c_{j})$ & $=$ & 
$\alpha _{M}^{P}(\sum_{j}m_{j}\otimes _{A}(1_{A}\otimes
_{S}c_{j})(\sum_{l}a_{l}\otimes g_{l}),$%
\end{tabular}%
\end{equation*}%
\emph{i.e.} $h\in \alpha _{M}^{P}(M\otimes _{A}(A\otimes _{S}C).$
Consequently, $\alpha _{M}^{P}$ is uniform.
\end{ex}

\begin{corollary}
\label{b-b}Let $B$ be an $S$-bisemialgebra.

\begin{enumerate}
\item If $_{S}B$ is an $\alpha $-semimodule, then we have isomorphisms of
categories%
\begin{equation*}
\mathbb{S}_{B}^{B}\simeq \mathbb{S}^{B\otimes _{S}^{a}B}\simeq \mathrm{Rat}%
^{B\otimes _{S}^{a}B}(\mathbb{S}_{B\#^{op}B^{\ast }})=\sigma \lbrack
B\otimes _{S}^{a}B_{B\#^{op}B^{\ast }}].
\end{equation*}

\item If $_{S}B$ is uniformly finitely presented and flat, then we have an
isomorphism of categories%
\begin{equation*}
\mathbb{S}_{B}^{B}\simeq \mathbb{S}_{B\#^{op}B^{\ast }}.
\end{equation*}
\end{enumerate}
\end{corollary}

\begin{Beweis}
\begin{enumerate}
\item This follows by Theorem \ref{HAC-iso} and Example \ref{ex-C-C}.

\item Since $_{S}B$ is uniformly finitely presented and flat, we have by
Lemma \ref{fp-flat} a canonical isomorphism $B\otimes _{S}B^{\ast }\overset{%
\upsilon _{(B,B,S)}}{\simeq }\mathrm{End}_{S}(B).$ An argument similar to
that of \cite[3.11]{BW2003} shows that $B^{\ast }$ is a left $B$%
-semicomodule through $B^{\ast }\overset{\chi }{\longrightarrow }\mathrm{End}%
_{S}(B)\overset{\upsilon _{(B,B,S)}}{\simeq }B\otimes _{S}B^{\ast },$ where $%
\chi (g)=\vartheta _{B}^{r}\circ (B\otimes _{S}g)\circ \Delta _{B}.$ Simple
computations show that the twisting map $\mathbf{\tau }_{(B,B)}:B^{\ast
}\otimes _{S}^{c}B\longrightarrow B\#^{op}B^{\ast }$ is right $%
B\#^{op}B^{\ast }$-linear, whence $B\#^{op}B^{\ast }\in \mathbb{S}_{B}^{B}$
and so $\mathbb{S}_{B}^{B}\simeq \mathbb{S}_{B\#^{op}B^{\ast }}.\blacksquare 
$
\end{enumerate}
\end{Beweis}

\section{Hopf Semialgebras}

\begin{punto}
With a \emph{Hopf }$S$\emph{-semialgebra}, we mean a datum $(H,\mu ,\eta
,\Delta ,\varepsilon ,\mathfrak{a}),$ where $(H,\mu ,\eta ,\Delta
,\varepsilon )$ is an $S$-bisemialgebra and $\mathrm{id}_{H}$ is a unit
(invertible) in the endomorphism semiring $(\mathrm{End}_{S}(H),\ast ),$ 
\emph{i.e. }there exists an $S$-linear map $\mathfrak{a}:H\longrightarrow $ $%
H$ such that%
\begin{equation*}
\sum \mathfrak{a}(h_{1})h_{2}=\varepsilon (h)1_{H}=\sum h_{1}\mathfrak{a}%
(h_{2})\text{ for all }h\in H.
\end{equation*}%
A \emph{morphism of Hopf }$S$\emph{-semialgebras} $f:H\longrightarrow
H^{\prime }$ is a morphism of $S$-bisemialgebras which is compatible with
the antipodes of $H$ and $H^{\prime }$ in the sense that%
\begin{equation}
\mathfrak{a}_{H^{\prime }}\circ f=f\circ \mathfrak{a}_{H}.  \label{Hopf-map}
\end{equation}%
The category of Hopf $S$-semialgebras is denoted by $\mathbf{HopfAlg}_{S}.$
\end{punto}

\begin{remarks}
\begin{enumerate}
\item If $(H,\mu ,\eta ,\Delta ,\varepsilon ,\mathfrak{a})$ is a Hopf $S$%
-semialgebra, then $\mathfrak{a}:H\longrightarrow H$ is a bialgebra
anti-morphism.

\item The category $\mathbf{HopAlg}_{S}\hookrightarrow \mathbf{Bialg}_{S}$
is a full subcategory, \emph{i.e.} if $H$ and $H^{\prime }$ are Hopf $S$%
-semialgebras and $f:H\longrightarrow H^{\prime }$ is an $S$-bialgebra
morphism, then $f$ is a morphism of Hopf $S$-semialgebras.

\item An $S$-bisemialgebra $B$ is a Hopf $S$-semialgebra with invertible
antipode $\mathfrak{a}$ if and only if $B^{cop}$ is a Hopf $S$-semialgebra
with invertible antipode $\widetilde{\mathfrak{a}}.$ Moreover, in this $%
\mathfrak{a}^{-1}=\widetilde{\mathfrak{a}}.$ In particular, if $H$ is a
commutative (cocommutative) Hopf $S$-semialgebra, then $\mathfrak{a}_{H}^{2}=%
\mathrm{id}$ (cf. \cite[Lemma 1.5.11, Corollary 1.5.12]{Mon1993}).
\end{enumerate}
\end{remarks}

\begin{definition}
Let $(H,\mu ,\eta ,\Delta ,\varepsilon ,\mathfrak{a})$ be a Hopf $S$%
-semialgebra. A bi-ideal $I\leq _{S}H$ is said to be a \emph{Hopf ideal} iff 
$\mathfrak{a}(I)\subseteq I.$
\end{definition}

\begin{proposition}
Let $(H,\mu ,\eta ,\Delta ,\varepsilon ,\mathfrak{a})$ be a Hopf $S$%
-semialgebra. For every Hopf ideal $I\leq _{S}H,$ we have a Hopf algebra
structure on $H/I;$ moreover, the canonical projection $\pi
_{I}:H\longrightarrow H/I$ is a morphism of Hopf $S$-semialgebras.
\end{proposition}

\begin{ex}
\label{G-group}Let $(G,\cdot ,e)$ be a group and consider the $S$%
-bisemialgebra $S[G]$ (Example \ref{S[M]}). One can easily see that%
\begin{equation*}
\mathfrak{a}:S[G]\longrightarrow S[G],\text{ }g\mapsto g^{-1}
\end{equation*}%
is an antipode for $S[G].$ So, $S[G]$ is a Hopf $S$-semialgebra.
\end{ex}

\begin{ex}
Consider the $S$-semialgebra $S[x,x^{-1}]$ with the usual multiplication and
unity. Notice that $G=\{x^{z}\mid z\in \mathbb{Z}\}$ is a group and that $%
S[x,x^{-1}]\simeq S[G]$ as $S$-semialgebras. It follows that $%
(S[x,x^{-1}],\mu ,\eta ,\Delta ,\varepsilon ,\mathfrak{a})$ is a Hopf $S$%
-semialgebra, where the structure maps are defined by extending the
following assignments linearly%
\begin{equation*}
\begin{tabular}{llllllll}
$\Delta $ & $:$ & $S[x,x^{-1}]$ & $\longrightarrow $ & $S[x,x^{-1}]\otimes
_{S}S[x,x^{-1}],$ & $x^{z}$ & $\mapsto $ & $x^{z}\otimes _{S}x^{z};$ \\ 
$\varepsilon $ & $:$ & $S[x,x^{-1}]$ & $\longrightarrow $ & $S,$ & $x^{z}$ & 
$\mapsto $ & $1_{S};$ \\ 
$\mathfrak{a}$ & $:$ & $S[x,x^{-1}]$ & $\longrightarrow $ & $S[x,x^{-1}],$ & 
$x^{z}$ & $\mapsto $ & $x^{-z}.$%
\end{tabular}%
\end{equation*}
\end{ex}

\begin{ex}
Let $(S,\oplus ,\cdot ,\mathbf{0},\mathbf{1})$ be an $S$-semiring which has
no non-zero zerodivisors and let $\mathbf{0}\neq \mathbf{2}=ab$ in $S,$ 
\emph{e.g.} $S$ is an additively idempotent semiring (Remark \ref{add-id}).
One can easily see that $H:=S[x]/(bx+x^{2})$ is a Hopf $S$-semialgebra with%
\begin{eqnarray*}
\Delta &:&H\longrightarrow H\otimes _{S}H,\text{ }\overline{x}\mapsto 
\overline{x}\otimes _{S}1_{H}\oplus 1_{H}\otimes _{S}\overline{x}\oplus a%
\overline{x}\otimes _{S}\overline{x}; \\
\varepsilon &:&H\longrightarrow S,\text{ }x\mapsto \mathbf{0}; \\
\mathfrak{a} &:&H\longrightarrow H,\text{ }\overline{x}\mapsto \overline{x}.
\end{eqnarray*}%
Notice that%
\begin{equation*}
\begin{tabular}{lllll}
$\sum \mathfrak{a}(\overline{x}_{1})\overline{x}_{2}$ & $=$ & $\mathfrak{a}(%
\overline{x})1_{H}\oplus \mathfrak{a}(1_{H})\overline{x}\oplus \mathfrak{a}(a%
\overline{x})\overline{x}$ & $=$ & $\overline{x}\oplus \overline{x}\oplus a%
\overline{x^{2}}$ \\ 
& $=$ & $\mathbf{2}\overline{x}\oplus a\overline{x^{2}}$ & $=$ & $(a\cdot b)%
\overline{x}\oplus a\overline{x^{2}}$ \\ 
& $=$ & $a(b\overline{x}\oplus \overline{x^{2}})$ & $=$ & $a(\overline{%
bx+x^{2}})$ \\ 
& $=$ & $a(0_{H})$ & $=$ & $0_{H}$ \\ 
& $=$ & $\mathbf{0(}1_{H})$ & $=$ & $\varepsilon (\overline{x})1_{H}.$%
\end{tabular}%
\end{equation*}%
Similarly, $\sum_{1}\overline{x}_{1}\mathfrak{a}(\overline{x}%
_{2})=\varepsilon (\overline{x})1_{H}.$ Consequently, $S[x]/(bx+x^{2})$ is a
Hopf $S$-semialgebra.$\blacksquare $
\end{ex}

\subsection*{Quantum Monoids}

\begin{definition}
A \emph{quantum monoid} is a non-commutative non-cocommutative Hopf
semialgebra.
\end{definition}

\begin{ex}
Consider the $S$-semialgebra with four \emph{different} generators%
\begin{equation*}
H=S<1,g,x,y\mid g^{2}=1,\text{ }x^{2}=xy=yx=y^{2}=0,\text{ }xg=gy,\text{ }%
yg=gx,\text{ }x+y=0>.
\end{equation*}%
Notice that $H$ is an $S$-bisemialgebra with comultiplication and counity
obtained by extending the following assignments as $S$-semialgebra morphisms%
\begin{eqnarray*}
\Delta (1) &=&1\otimes _{S}1,\text{ }\Delta (g)=g\otimes _{S}g,\text{ }%
\Delta (x)=x\otimes _{S}1+g\otimes _{S}x,\text{ }\Delta (y)=y\otimes
_{S}1+g\otimes _{S}y; \\
\varepsilon (1) &=&1=\varepsilon (g),\text{ }\varepsilon (x)=0=\varepsilon
(y).
\end{eqnarray*}%
Moreover, $H$ is a Hopf $S$-semialgebra with antipode defined by extending
the following assignments linearly%
\begin{equation*}
\mathfrak{a}(1)=1,\text{ }\mathfrak{a}(g)=g,\text{ }\mathfrak{a}(x)=xg,\text{
}\mathfrak{a}(y)=yg.
\end{equation*}%
Clearly, $H$ is \emph{non-commutative} and \emph{non-cocommutative}, i.e. $H$
is a quantum monoids. Notice that $H$ is in fact a semialgebraic version of
Sweedler's Hopf Algebra (quantum group).
\end{ex}

\begin{ex}
Let $n\geq 2$ and $q\in S$ be such that $q^{n}=1_{S}$ and $q^{i}\neq 1$ for
any $i\in \{1,\cdots ,n-1\}.$ Consider the $S$-semialgebra with four \emph{%
different} generators%
\begin{equation*}
H=S<1,g,x,y\mid g^{n}=1,\text{ }x^{i}y^{n-i}=0=y^{n-i}x^{i},\text{ }%
i=0,\cdots ,n,\text{ }xg=qgx,\text{ }yg=qgy,\text{ }x+y=0>.
\end{equation*}%
Notice that $H$ is an $S$-bisemialgebra with comultiplication and counity
obtained by extending the following assignments as $S$-semialgebra morphisms%
\begin{eqnarray*}
\Delta (1) &=&1\otimes _{S}1,\text{ }\Delta (g)=g\otimes _{S}g,\text{ }%
\Delta (x)=x\otimes _{S}1+g\otimes _{S}x,\text{ }\Delta (y)=y\otimes
_{S}1+g\otimes _{S}y; \\
\varepsilon (1) &=&1=\varepsilon (g),\text{ }\varepsilon (x)=0=\varepsilon
(y).
\end{eqnarray*}%
Moreover, $H$ is a Hopf $S$-semialgebra with antipode defined by extending
the following assignments linearly%
\begin{equation*}
\mathfrak{a}(1)=1,\text{ }\mathfrak{a}(g)=g^{-1},\text{ }\mathfrak{a}%
(x)=g^{-1}y,\text{ }\mathfrak{a}(y)=g^{-1}x.
\end{equation*}%
Clearly, $H$ is \emph{non-commutative} and \emph{non-cocommutative}, i.e. $H$
is a quantum monoid. Notice that $H$ is in fact a semialgebraic version of
Taft's Hopf Algebra (quantum group).
\end{ex}

\begin{ex}
The $S$-semialgebra $H=S[x,y,y^{-1}]/(xy+yx,x^{2}),$ where $x$ and $y$ are 
\emph{non-commuting} indeterminates, is a quantum monoid with%
\begin{eqnarray*}
\Delta (1) &=&1\otimes _{S}1,\text{ }\Delta (x)=x\otimes _{S}1+y^{-1}\otimes
_{S}x,\text{ }\Delta (y)=y\otimes _{S}y; \\
\varepsilon (1) &=&1,\text{ }\varepsilon (x)=0,\text{ }\varepsilon (y)=1; \\
\mathfrak{a}(1) &=&1,\text{ }\mathfrak{a}(x)=xy,\text{ }\mathfrak{a}%
(y)=y^{-1},\text{ }\mathfrak{a}(y^{-1})=y.
\end{eqnarray*}%
In fact, $H$ is a semialgebraic version of Pareigis' Hopf algebra (quantum
group) \cite[Example 3.4.22]{HGK2010}.
\end{ex}

\subsection*{Fundamental Theorem}

In what follows we give sufficient and necessary conditions for a given $S$%
-bisemialgebra to have an antipode (\emph{i.e.} to be a Hopf $S$%
-semialgebra).

\begin{proposition}
\label{antipode}\-Let $B$ be an $S$-bisemialgebra. The following are
equivalent:

\begin{enumerate}
\item $B$ is a Hopf $S$-semialgebra;

\item $B\otimes _{S}B\overset{\gamma _{B}}{\simeq }B\otimes _{S}^{a}B$ is an
isomorphism in $\mathbb{S}_{B}^{B},$ where $\gamma _{B}(a\otimes _{S}b)=\sum
ab_{1}\otimes _{S}^{a}b_{2};$

\item $B\otimes _{S}^{c}B\overset{\gamma ^{B}}{\simeq }B\otimes _{S}B$ is an
isomorphism in $\mathbb{S}_{B}^{B},$ where $\gamma ^{B}(a\otimes
_{S}^{c}b)=\sum a_{1}\otimes _{S}a_{2}b.$
\end{enumerate}
\end{proposition}

\begin{Beweis}
The proof is based on direct standard computations (\emph{e.g. }\cite%
{Ion1998}). Since, our semimodules do not allow subtraction (in general), we
notice here that the proof of \cite[15.2 (3)]{BW2003} does \emph{not} work
in our case. We prove $(1)\Longleftrightarrow (2);$ notice that $%
(1)\Longleftrightarrow (3)$ can be proved symmetrically.

$(1)\Rightarrow (2)$ If $\mathfrak{a}$ is an antipode for $B,$ then it is
clear that $\gamma _{B}$ is an isomorphism with inverse%
\begin{equation*}
\omega _{B}:B\otimes _{S}^{a}B\longrightarrow B\otimes _{S}B,\text{ }%
a\otimes _{S}^{a}b\mapsto \sum a\mathfrak{a}(b_{1})\otimes _{S}b_{2}.
\end{equation*}

$(2)\Rightarrow (1)$ Assume that $\gamma _{B}$ is invertible with inverse $%
\omega _{B}:B\otimes _{S}B\longrightarrow B\otimes _{S}^{a}B$ in $\mathbb{S}%
_{B}^{B};$ in particular, $\omega _{B}$ is right $B$-colinear, \emph{i.e.}%
\begin{equation}
(\omega _{B}\otimes _{S}B)\circ (B\otimes _{S}\Delta )=(B\otimes _{S}\Delta
)\circ \omega _{B}  \label{B-colinear}
\end{equation}%
Moreover, since $\gamma _{B}$ is left $B$-linear, its inverse $\omega _{B}$
is indeed left $B$-linear. Setting%
\begin{equation*}
\xi _{B}^{r}:B\otimes _{S}B\overset{B\otimes _{S}\varepsilon }{%
\longrightarrow }B\otimes _{S}S\overset{\vartheta _{B}^{r}}{\simeq }B\text{
and }\omega _{B}(1\otimes _{S}b)=\sum b^{(0)}\otimes _{S}b^{(1)},
\end{equation*}%
we claim that%
\begin{equation*}
\mathfrak{a}:B\longrightarrow B,\text{ }b\mapsto (\xi _{B}^{r}\circ \omega
_{B})(1\otimes _{S}b)=\sum b^{(1)}\varepsilon (b^{(2)})
\end{equation*}%
is an antipode for $B.$

On one hand, we have%
\begin{equation*}
\begin{tabular}{lll}
&  & $(\mu \circ \mathfrak{a}\otimes _{S}B\circ \Delta )(b)$ \\ 
& $=$ & $\sum \mathfrak{a}(b_{1})b_{2}$ \\ 
& $=$ & $\sum \mathfrak{a}(b_{1})b_{2}\varepsilon (b_{3})$ \\ 
& $=$ & $[\mu \circ B\otimes _{S}\xi _{B}^{r}](\sum \mathfrak{a}%
(b_{1})\otimes _{S}b_{2}\otimes _{S}b_{3})$ \\ 
& $=$ & $[\mu \circ B\otimes _{S}\xi _{B}^{r}](\sum (b_{1})^{(0)}\varepsilon
((b_{1})^{(1)})\otimes _{S}b_{2}\otimes _{S}b_{3})$ \\ 
& $=$ & $[\mu \circ B\otimes _{S}\xi _{B}^{r}\circ \xi _{B}^{r}\otimes
_{S}B\otimes _{S}B](\sum (b_{1})^{(0)}\otimes _{S}(b_{1})^{(1)}\otimes
_{S}b_{2}\otimes _{S}b_{3})$ \\ 
& $=$ & $[\mu \circ B\otimes _{S}\xi _{B}^{r}\circ \xi _{B}^{r}\otimes
_{S}B\otimes _{S}B\circ B\otimes _{S}B\otimes _{S}\Delta ](\sum
(b_{1})^{(0)}\otimes _{S}(b_{1})^{(1)}\otimes _{S}b_{2})$ \\ 
& $=$ & $[\mu \circ B\otimes _{S}\xi _{B}^{r}\circ \xi _{B}^{r}\otimes
_{S}B\otimes _{S}B\circ B\otimes _{S}B\otimes _{S}\Delta \circ \omega
_{B}\otimes _{S}B](\sum (1\otimes _{S}b_{1})\otimes _{S}b_{2})$ \\ 
& $=$ & $[\mu \circ B\otimes _{S}\xi _{B}^{r}\circ \xi _{B}^{r}\otimes
_{S}B\otimes _{S}B\circ B\otimes _{S}B\otimes _{S}\Delta \circ \omega
_{B}\otimes _{S}B\circ B\otimes _{S}\Delta ](1\otimes _{S}b)$ \\ 
& $\overset{(\text{\ref{B-colinear}})}{=}$ & $[\mu \circ B\otimes _{S}\xi
_{B}^{r}\circ \xi _{B}^{r}\otimes _{S}B\otimes _{S}B\circ B\otimes
_{S}B\otimes _{S}\Delta \circ B\otimes _{S}\Delta \circ \omega
_{B}](1\otimes _{S}b)$ \\ 
& $=$ & $[\mu \circ B\otimes _{S}\xi _{B}^{r}\circ \xi _{B}^{r}\otimes
_{S}B\otimes _{S}B\circ B\otimes _{S}B\otimes _{S}\Delta \circ B\otimes
_{S}\Delta ](\sum b^{(0)}\otimes _{S}b^{(1)})$ \\ 
& $=$ & $[\mu \circ B\otimes _{S}\xi _{B}^{r}\circ \xi _{B}^{r}\otimes
_{S}B\otimes _{S}B\circ B\otimes _{S}B\otimes _{S}\Delta ](\sum
b^{(0)}\otimes _{S}(b^{(1)})_{1}\otimes _{S}(b^{(1)})_{2})$ \\ 
& $=$ & $[\mu \circ B\otimes _{S}\xi _{B}^{r}\circ \xi _{B}^{r}\otimes
_{S}B\otimes _{S}B](\sum b^{(0)}\otimes _{S}(b^{(1)})_{1}\otimes
_{S}(b^{(1)})_{2}\otimes _{S}(b^{(1)})_{3}$ \\ 
& $=$ & $[\mu \circ B\otimes _{S}\xi _{B}^{r}](\sum b^{(0)}\varepsilon
((b^{(1)})_{1})\otimes _{S}(b^{(1)})_{2}\otimes _{S}(b^{(1)})_{3}$ \\ 
& $=$ & $[\mu \circ B\otimes _{S}\xi _{B}^{r}](\sum b^{(0)}\otimes
_{S}\varepsilon ((b^{(1)})_{1})(b^{(1)})_{2}\otimes _{S}(b^{(1)})_{3})$ \\ 
& $=$ & $[\mu \circ B\otimes _{S}\xi _{B}^{r}](\sum b^{(0)}\otimes
_{S}(b^{(1)})_{1}\otimes _{S}(b^{(1)})_{2})$ \\ 
& $=$ & $\mu (\sum b^{(0)}\otimes _{S}(b^{(1)})_{1}\varepsilon
((b^{(1)})_{2}))$ \\ 
& $=$ & $\sum b^{(0)}(b^{(1)})_{1}\varepsilon ((b^{(1)})_{2})$ \\ 
& $=$ & $\xi _{B}^{r}(\sum b^{(0)}(b^{(1)})_{1}\otimes _{S}(b^{(1)})_{2})$
\\ 
& $=$ & $(\xi _{B}^{r}\circ \gamma _{B})(\sum b^{(0)}\otimes _{S}b^{(1)})$
\\ 
& $=$ & $(\xi _{B}^{r}\circ \gamma _{B}\circ \omega _{B})(1\otimes _{S}b)$
\\ 
& $=$ & $\xi _{B}^{r}(1\otimes _{S}b)$ \\ 
& $=$ & $1_{B}\varepsilon (b)$ \\ 
& $=$ & $(\eta \circ \varepsilon )(b)$%
\end{tabular}%
\end{equation*}%
On the other hand, we have%
\begin{equation*}
\begin{tabular}{lllll}
$(\mu \circ B\otimes _{S}\mathfrak{a}\circ \Delta )(b)$ & $=$ & $\sum b_{1}%
\mathfrak{a}(b_{2})$ &  &  \\ 
& $=$ & $\sum b_{1}(b_{2})^{(0)}\varepsilon ((b_{2})^{(1)})$ &  &  \\ 
& $=$ & $\xi _{B}^{r}(\sum b_{1}(b_{2})^{(0)}\otimes _{S}(b_{2})^{(1)})$ & 
&  \\ 
& $=$ & $\xi _{B}^{r}(\sum b_{1}\omega _{B}(1\otimes _{S}b_{2})$ &  &  \\ 
& $=$ & $[\xi _{B}^{r}\circ \omega _{B}](\sum b_{1}\otimes _{S}b_{2})$ &  & (%
$\omega _{B}$ is left $B$-linear) \\ 
& $=$ & $[\xi _{B}^{r}\circ \omega _{B}\circ \gamma _{B}](1\otimes _{S}b)$ & 
&  \\ 
& $=$ & $\xi _{B}^{r}(1\otimes _{S}b)$ &  & ($\omega _{B}\circ \gamma _{B}=%
\mathrm{id}_{B\otimes _{S}B}$) \\ 
& $=$ & $1_{B}\varepsilon (b)$ &  &  \\ 
& $=$ & $(\eta \circ \varepsilon )(b).\blacksquare $ &  & 
\end{tabular}%
\end{equation*}
\end{Beweis}

We present now the \emph{Fundamental Theorem of Hopf Semimodules}.

\begin{theorem}
\label{FTHM}The following are equivalent for an $S$-bisemialgebra $(B,\mu
,\eta ,\Delta ,\varepsilon ):$

\begin{enumerate}
\item $B$ is a Hopf $S$-semialgebra;

\item For every $M\in \mathbb{S}_{B}^{B},$ we have an isomorphism in $%
\mathbb{S}_{B}^{B}:$%
\begin{equation}
M^{\mathrm{co}B}\otimes _{S}B\overset{\psi _{M}}{\simeq }M,\text{ }m\otimes
_{S}b\mapsto mb;  \label{co-id}
\end{equation}

\item For every $M\in \mathbb{S}_{B}^{B},$ we have an isomorphism in $%
\mathbb{S}_{B}^{B}:$%
\begin{equation}
\mathrm{Hom}_{B}^{B}(B,M)\otimes _{S}B\overset{\varphi _{M}}{\simeq }M,\text{
}f\otimes _{S}b\mapsto f(b);  \label{BB-id}
\end{equation}

\item We have an isomorphism in $\mathbb{S}_{B}^{B}:$%
\begin{equation}
\mathrm{Hom}_{B}^{B}(B,B\otimes _{S}^{a}B)\otimes _{S}B\overset{\varphi
_{B\otimes _{S}B}}{\simeq }B\otimes _{S}^{a}B.  \label{BB}
\end{equation}
\end{enumerate}
\end{theorem}

\begin{Beweis}
$(1)\Rightarrow (2):$ Let $\mathfrak{a}$ be an antipode for $H.$ It can be
shown that $\psi _{M}$ is an isomorphism with inverse given by $m\mapsto
\sum m_{<0>}\mathfrak{a}(m_{<1>}).$

$(2)\Rightarrow (3)$ This follows from the isomorphisms%
\begin{equation*}
\mathrm{Hom}_{B}^{B}(B,M)\otimes _{S}B\overset{\upsilon _{M}\otimes _{S}B}{%
\underset{(\ref{f(1)})}{\simeq }}M^{\mathrm{co}B}\otimes _{S}B\overset{\psi
_{M}}{\simeq }M
\end{equation*}%
and the fact that $\psi _{M}\circ \upsilon _{M}\otimes _{S}B=\varphi _{M}:$
for all $f\in \mathrm{Hom}_{B}^{B}(B,M)$ and $b\in B:$%
\begin{equation*}
(\psi _{M}\circ (\upsilon _{M}\otimes _{S}B))(f\otimes _{S}b)=\psi
_{M}(f(1_{B})\otimes _{S}b)=f(1_{B})b=f(b)=\varphi _{M}(f\otimes _{S}b).
\end{equation*}

$(3)\Rightarrow (4)$ trivial.

$(4)\Rightarrow (1)$ We have an isomorphism of $S$-semimodules%
\begin{equation*}
B\overset{\zeta _{B}}{\simeq }\mathrm{Hom}_{B}^{B}(B,B\otimes _{S}^{a}B),%
\text{ }a\mapsto \lbrack b\mapsto \sum ab_{1}\otimes _{S}^{a}b_{2}]
\end{equation*}%
with inverse $f\mapsto (\vartheta _{B}^{r}\circ B\otimes _{S}\varepsilon
\circ f)(1_{B}).$ Moreover, for all $a,b\in B$ we have%
\begin{equation*}
(\varphi _{B\otimes _{S}B}\circ \zeta _{B}\otimes _{S}B)(a\otimes
_{S}b)=\varphi _{B\otimes _{S}B}(\zeta _{B}(a)\otimes _{S}b)=\zeta
_{B}(a)(b)=\sum ab_{1}\otimes _{S}^{a}b_{2}=\gamma _{B}(a\otimes _{S}b).
\end{equation*}%
Consequently, we have an isomorphism in $\mathbb{S}_{B}^{B}:$%
\begin{equation*}
\gamma _{B}:B\otimes _{S}B\overset{\zeta _{B}\otimes _{S}B}{\simeq }\mathrm{%
Hom}_{B}^{B}(B,B\otimes _{S}^{a}B)\otimes _{S}B\overset{\varphi _{B\otimes
_{S}B}}{\simeq }B\otimes _{S}^{a}B.
\end{equation*}%
By Proposition \ref{antipode}, $H$ is a Hopf $S$-semialgebra.$\blacksquare $
\end{Beweis}

We are ready now to present the \emph{Fundamental Theorem of Hopf
Semialgebras}. Making use of Proposition \ref{antipode}, the proof is
similar to that of \cite[15.5]{BW2003}.

\begin{theorem}
\label{FTHA}The following are equivalent for an $S$-bisemialgebra $B:$

\begin{enumerate}
\item $B$ is a Hopf $S$-semialgebra;

\item $\mathbb{S}_{S}\overset{-\otimes _{S}B}{\approx }\mathbb{S}_{B}^{B}$ 
\emph{(}an equivalence of categories\emph{)} with inverse $\mathrm{Hom}%
_{B}^{B}(B,-):\mathbb{S}_{B}^{B}\longrightarrow \mathbb{S}_{S}.$
\end{enumerate}
\end{theorem}

\begin{Beweis}
$(1)\Rightarrow (2)$ For every $M_{S},$ it is obvious that we have a natural
isomorphism of $S$-semimodules%
\begin{equation*}
M\simeq \mathrm{Hom}_{B}^{B}(B,M\otimes _{S}B),\text{ }m\mapsto \lbrack
b\mapsto m\otimes _{S}b]
\end{equation*}%
with inverse $f\mapsto (\rho _{M}\circ f)(1_{B}).$ It follows that $\mathrm{%
Hom}_{B}^{B}(B,-)\circ -\otimes _{S}B\simeq \mathrm{id}_{\mathbb{S}_{S}}.$
On the other hand, we have a natural isomorphism $\mathrm{Hom}%
_{B}^{B}(B,M)\otimes _{S}B\overset{\varphi _{M}}{\simeq }M$ in $\mathbb{S}%
_{B}^{B};$ this means that $-\otimes _{S}B\circ \mathrm{Hom}%
_{B}^{B}(B,-)\simeq \mathrm{id}_{\mathbb{S}_{B}^{B}}.$ Consequently, $%
\mathbb{S}_{B}^{B}\overset{\mathrm{Hom}_{B}^{B}(B,-)}{\simeq }\mathbb{S}_{S}$
with inverse $-\otimes _{S}B.$

$(2)\Rightarrow (1)$ This follows by Proposition \ref{FTHM}.$\blacksquare $
\end{Beweis}

\begin{remark}
We notice here that the Fundamental Theorem of Hopf Algebras can be obtained
by applying results of \cite{BW2011} on \emph{Hopf monads} taking into
consideration that $\mathbb{S}_{S}$ has colimits and that $-\otimes
_{S}B\simeq B\otimes _{S}-:\mathbb{S}_{S}\longrightarrow \mathbb{S}_{S}$
preserves colimits for every $B_{S}.$ However, we follow in this paper the
direct algebraic approach.
\end{remark}

We present now the main reconstruction result in this Section:

\begin{theorem}
\label{H-monad}Let $H$ be an $S$-semimodule. There is a bijective
correspondence between the structures of Hopf $S$-semialgebras on $H,$ the
Hopf monad structures on $H\otimes _{S}-:\mathbb{S}_{S}\longrightarrow 
\mathbb{S}_{S}$ and the Hopf monad structures on $-\otimes _{S}H:\,\mathbb{S}%
_{S}\longrightarrow \,\mathbb{S}_{S}.$
\end{theorem}

\begin{Beweis}
The proof of the bijective correspondence in (2) is similar to that of \cite[%
Theorem 3.9]{Ver} taking into consideration that $\mathbb{S}_{S}$ is
cocomplete, that $S$ is a regular generator in $\mathbb{S}_{S}$ and the fact
that $X\otimes _{S}-\simeq -\otimes _{S}X$ preserve colimits in $\mathbb{S}%
_{S}$ for every $S$-semimodule $X.\blacksquare $
\end{Beweis}

\subsection*{Semisimple and Cosemisimple Hopf Semialgebras}

\begin{definition}
Let $\mathfrak{C}$ be a category. We say that an object $V$ in $\mathfrak{C}$
is \emph{semisimple} (\emph{simple}) iff every monomorphism $\iota
:U\longrightarrow V$ in $\mathfrak{C}$ is a coretraction (an isomorphism).
Moreover, we say that $V$ is \emph{completely reducible} iff $V$ is a direct
sum of simple objects in $\mathfrak{C}.$
\end{definition}

\begin{remark}
In contrast with modules over a ring, a semisimple semimodule over a
semiring is not necessarily completely reducible:\ an $S$-subsemimodule $U%
\overset{\iota }{\hookrightarrow }V$ for which $\iota $ is a coretraction is
not necessarily a summand \cite{Gol1999}.
\end{remark}

\begin{definition}
Let $\mathfrak{C}$ and $\mathfrak{D}$ be two categories such that $\mathrm{%
Obj}(\mathfrak{C})\subseteq \mathrm{Obj}(\mathfrak{D})$ and $\mathrm{Mor}_{%
\mathfrak{C}}(X,Y)\subseteq \mathrm{Mor}_{\mathfrak{D}}(X,Y)$ for all $%
X,Y\in \mathrm{Obj}(\mathfrak{C}).$ We say that $V\in \mathfrak{C}$ is $(%
\mathfrak{C},\mathfrak{D})$\emph{-semisimple} iff every monomorphism $%
U\longrightarrow V$ in $\mathfrak{C}$ which is a coretraction in $\mathfrak{D%
}$ is also a coretraction in $\mathfrak{C}.$
\end{definition}

\begin{definition}
Let $H$ b e a Hopf $S$-semialgebra.

\begin{enumerate}
\item We say that a left $H$-semimodule $M$ is

\emph{semisimple} iff $M$ is semisimple in $_{H}\mathbb{S};$

$(H,S)$\emph{-semisimple} iff $M$ is $(_{H}\mathbb{S},\mathbb{S}_{S})$%
-semisimple.

\item We say that a left $H$-semicomodule $N$ is

\emph{cosemisimple} iff $N$ is semisimple in $^{H}\mathbb{S};$

$(H,S)$\emph{-cosemisimple} iff $N$ is $(^{H}\mathbb{S},\mathbb{S}_{S})$%
-semisimple.
\end{enumerate}
\end{definition}

\begin{definition}
We say that a Hopf $S$-semialgebra $H$ is

\emph{left semisimple} ($(H,S)$-\emph{semisimple}) iff $H\in $ $_{H}\mathbb{S%
}$ is semisimple ($(H,S)$-semisimple);

\emph{left cosemisimple} ($(H,S)$-\emph{cosemisimple}) iff $H\in $ $^{H}%
\mathbb{S}$ is cosemisimple ($(H,S)$-cosemisimple);

\emph{separable} iff $H$ is a separable $S$-semialgebra (\emph{i.e.} iff $%
\mu _{H}:H\otimes _{S}H\longrightarrow H$ splits in $_{H}S_{H}$);

\emph{coseparable} iff $H$ is a coseparable $S$-semicoalgebra (\emph{i.e.}
iff $\Delta _{H}:H\longrightarrow H\otimes _{S}H$ splits in $^{H}S^{H}$).

Symmetrically, one can introduce right-sided versions of these notions.
\end{definition}

\begin{proposition}
\label{int-imp}Let $(H,\mu ,\eta ,\Delta ,\varepsilon ,a)$ be a Hopf $S$%
-semialgebra with $_{S}H$ an $\alpha $-semimodule and consider $H^{\ast
rat}:=\mathrm{Rat}^{H}(_{H^{\ast }}H^{\ast }).$

\begin{enumerate}
\item $B^{\ast rat}\in \mathbb{S}_{H}^{H}$ with right $H$-action given by $%
(f\leftharpoondown h)(g):=f(\mathfrak{a}(h)g)$ for all $f\in B^{\ast rat}$
and $h,g\in H.$

\item We have an isomorphism in $\mathbb{S}_{H}^{H}:$%
\begin{equation*}
\int^{l}H\otimes _{S}H\longrightarrow H^{\ast rat},\text{ }t\otimes
_{S}h\mapsto t\leftharpoondown h.
\end{equation*}

\item $\int_{l}H=0$ if and only if $H^{\ast rat}=0.$
\end{enumerate}
\end{proposition}

\begin{Beweis}
\begin{enumerate}
\item Direct calculations \cite[5.2.1]{DNR2001} show that $H^{\ast rat}\in 
\mathbb{S}_{H}^{H}.$

\item This follows directly from the fact that $\int^{l}H=(H^{\ast rat})^{%
\mathrm{co}H}$ (Lemma \ref{left-int}) and the Fundamental Theorem of Hopf
Semimodules (Proposition \ref{FTHM}) applied to $M=H^{\ast rat}.$

\item If $H^{\ast rat}=0,$ then $\int^{l}H=(H^{\ast rat})^{\mathrm{co}H}=0.$
On the other hand, if $\int^{l}H=0,$ then $H^{\ast rat}\simeq
\int^{l}H\otimes _{S}H=0\otimes _{S}H=0.\blacksquare $
\end{enumerate}
\end{Beweis}

The following result characterizes coseparable Hopf $S$-semialgebras.

\begin{proposition}
\label{H-cosemi}The following are equivalent for a Hopf $S$-semialgebra $H:$

\begin{enumerate}
\item $\int^{l,1}H\neq \varnothing ;$

\item Every left \emph{(}right\emph{)} $H$-semicomodule is $(H,S)$%
-cosemisimple;

\item $H$ is left \emph{(}right\emph{)} $(H,S)$-cosemisimple;

\item $H$ is coseparable;

\item There exists an $S$-linear map $\delta :H\otimes _{S}H\longrightarrow
S $ such that%
\begin{equation*}
\delta \circ \Delta _{H}=\varepsilon \text{ and }(H\otimes _{S}\delta )\circ
(\Delta _{H}\otimes _{S}H)=(\delta \otimes _{S}H)\circ (H\otimes _{S}\Delta
_{H});
\end{equation*}

\item $\int^{r,1}H\neq \varnothing .$
\end{enumerate}
\end{proposition}

\begin{Beweis}
The proof can be obtained by arguments similar to a combination of \cite[%
16.10]{BW2003} and \cite[16.10]{BW2003} (see also \cite[Theorem 3.3.2]%
{Abe1980}).$\blacksquare $
\end{Beweis}

\begin{corollary}
\label{cos-rat}Let $H$ be a Hopf $S$-semialgebra. If $_{S}H$ is an $\alpha $%
-semimodule, then the following are equivalent:

\begin{enumerate}
\item $H$ is left $(H,S)$-cosemisimple;

\item $H$ is a semisimple right $H^{\ast }$-semimodule;

\item every rational right $H^{\ast }$-semimodule is semisimple.
\end{enumerate}
\end{corollary}

\begin{Beweis}
If $_{S}H$ is an $\alpha $-semimodule, then we have an isomorphism of
categories $^{H}\mathbb{S}\simeq $ $^{H}\mathrm{Rat}(\mathbb{S}_{H^{\ast }})$
and the result follows by Proposition \ref{H-cosemi}.$\blacksquare $
\end{Beweis}

\begin{ex}
Let $S$ be such that every $S$-subsemimodule of $H$ is injective. If $H$ has
a total left integral, then $H$ is left cosemisimple by Proposition \ref%
{H-cosemi}.
\end{ex}

The following result characterizes separable Hopf $S$-semialgebras.

\begin{proposition}
\label{H-semi}The following are equivalent for a Hopf $S$-semialgebra $H:$

\begin{enumerate}
\item $\int_{l,1}H\neq \varnothing ;$

\item Every left \emph{(}right\emph{)} $H$-semimodule is $(H,S)$-semisimple;

\item $H$ is left \emph{(}right\emph{)} $(H,S)$-semisimple;

\item $H$ is separable;

\item There exists $e=\sum e^{1}\otimes _{S}e^{2}\in H\otimes _{S}H$ such
that%
\begin{equation*}
he=eh\text{ for all }h\in H\text{ and }\sum e^{1}e^{2}=1_{H};
\end{equation*}

\item $\int_{r,1}H\neq \varnothing .$
\end{enumerate}
\end{proposition}

\begin{Beweis}
The proof is similar to that of \cite[16.13]{BW2003} (see also \cite{CMZ2002}%
).$\blacksquare $
\end{Beweis}

\section{Dual Bisemialgebras and Dual Hopf}

\qquad As before, $S$ denotes a commutative semiring with $1_{S}\neq 0_{S}$
and $\mathbb{S}_{S}$ is the category of $S$-semimodules.

\begin{lemma}
\label{X*}Let $X,Y$ be $S$-semimodules. If $X_{S}$ is finitely generated and
projective and $Y_{S}$ is uniformly finitely presented, then we have a
canonical isomorphism%
\begin{equation}
X^{\ast }\otimes _{S}Y^{\ast }\overset{\mathfrak{h}_{(X,Y)}}{\simeq }%
(X\otimes _{S}X)^{\ast },\text{ }f\otimes _{S}g\mapsto \lbrack x\otimes
_{S}y\mapsto f(x)g(y)].  \label{*x*=*}
\end{equation}
\end{lemma}

\begin{Beweis}
Since $X_{S}$ is finitely generated and projective, $X_{S}^{\ast }$ is also
(finitely generated and) projective whence flat. The given map is the
composition of the following isomorphisms%
\begin{equation*}
X^{\ast }\otimes _{S}\mathrm{Hom}_{S}(Y,S)\overset{\upsilon _{(X^{\ast
},Y,S)}}{\simeq }\mathrm{Hom}_{S}(X,\mathrm{Hom}_{S}(Y,S)\otimes
_{S}S)\simeq \mathrm{Hom}_{S}(X\otimes _{S}Y,S).
\end{equation*}%
where $\upsilon _{(X^{\ast },Y,S)}$ is an isomorphism by Lemma \ref{fp-flat}
and the second isomorphism is the canonical one indicating the adjointness
of the tensor and hom functors.$\blacksquare $
\end{Beweis}

\begin{proposition}
\label{A*-co}Let $A$ be an $S$-semimodule with $A_{S}$ uniformly finitely
presented projective and consider the isomorphism of semirings $S^{\ast }%
\overset{h}{\simeq }S.$

\begin{enumerate}
\item If $(A,\mu ,\eta )$ is an $S$-semialgebra, then $(A^{\ast },\mathfrak{h%
}_{(A,A)}^{-1}\circ \mu ^{\ast },h\circ \eta ^{\ast })$ is an $S$%
-semicoalgebra.

\item If $(A,\mu ,\eta ,\Delta ,\varepsilon )$ is an $S$-bisemialgebra, then 
$(A^{\ast },\Delta ^{\ast }\circ \mathfrak{h}_{(A,A)},\varepsilon ^{\ast
}\circ h^{-1},\mathfrak{h}_{(A,A)}^{-1}\circ \mu ^{\ast },h\circ \eta ^{\ast
}),$ is an $S$-bisemialgebra.

\item If $(A,\mu ,\eta ,\Delta ,\varepsilon ,\mathfrak{a})$ is a Hopf $S$%
-semialgebra, then $(A^{\ast },\Delta ^{\ast }\circ \mathfrak{h}%
_{(A,A)},\varepsilon ^{\ast }\circ h^{-1},\mathfrak{h}_{(A,A)}^{-1}\circ \mu
^{\ast },h\circ \eta ^{\ast },\mathfrak{a}^{\ast })$ is a Hopf $S$%
-semialgebra.
\end{enumerate}
\end{proposition}

\begin{Beweis}
Applying Lemma \ref{X*} twice, we have the canonical isomorphisms%
\begin{equation*}
(A\otimes _{S}A)^{\ast }\overset{\mathfrak{h}_{(A,A)}}{\simeq }A^{\ast
}\otimes _{S}A^{\ast }\text{ and }(A\otimes _{S}A\otimes _{S}A)^{\ast }%
\overset{\mathfrak{h}_{(A\otimes _{S}A,A)}}{\simeq }(A\otimes _{S}A)^{\ast
}\otimes _{S}A^{\ast }\overset{\mathfrak{h}_{(A,A)}\otimes _{S}A^{\ast }}{%
\simeq }A^{\ast }\otimes _{S}A^{\ast }\otimes _{S}A^{\ast }.
\end{equation*}%
The result follows now by applying $(-)^{\ast }:=\mathrm{Hom}_{S}(-,S)$ to
the $S$-semimodules and arrows defining the given $S$-semialgebra \emph{(}%
resp. $S$-bisemialgebra, Hopf $S$-semialgebra).
\end{Beweis}

\subsection*{The finite dual}

\begin{definition}
We say that an $S$-semimodule $M$ is Noetherian iff every $S$-subsemimodule
of $M$ is finitely generated. We say that the semiring $S$ is right (left)
Noetherian iff $S_{S}$ ($_{S}S$) is Noetherian. If $S$ is right and left
Noetherian, then we say that $S$ is Noetherian.
\end{definition}

\begin{definition}
We say that the semiring $S$ is \emph{right }(\emph{left}) $m$-\emph{%
Noetherian} iff every finitely generated right (left) $S$-semimodule is
Noetherian. If $S$ is right and left $m$-Noetherian, then we say that $S$ is
an $m$\emph{-Noetherian semiring}.
\end{definition}

\begin{remark}
The notion of $m$-Noetherian semirings was considered in \cite{EM2011} where
such semirings were called Noetherian. To avoid confusion, we add the prefix 
$m$ which refers to modules.
\end{remark}

\begin{ex}
Every \emph{finite} semiring is $m$-Noetherian. Such rings are not rare: Let 
$n$ be a positive integer and $X_{n}=\{-\infty ,0,1,\cdots ,n\}.$ Define
addition and multiplication on $X_{n}$ as%
\begin{equation*}
i+_{X_{n}}h:=\max \{i,h\}\text{ and }i\cdot _{X_{n}}h=\min \{i+h,n\}.
\end{equation*}%
One can easily check that $\{X_{n},+_{X_{n}},\cdot _{X_{n}},\mathbf{0},%
\mathbf{1}\}$ is a (finite) semiring with $\mathbf{0}=-\infty $ and $\mathbf{%
1}=0.$ Moreover, $X_{n}$ has no non-zero zerodivisors and is additively
idempotent since%
\begin{equation*}
\mathbf{1}+_{X_{n}}\mathbf{1}=0+_{X_{n}}0=\max \{0,0\}=0=\mathbf{1}.
\end{equation*}
\end{ex}

\begin{ex}
(\cite[Example 2.2]{BMS2013})\ The semiring $(\mathbb{N}_{0},+,\cdot )$ is
Noetherian but not $m$-Noetherian: the semimodule $\mathbb{N}_{0}\times 
\mathbb{N}_{0}$ is finitely generated but its subsemimodule generated by $%
\{(n+1,n)\mid n\in \mathbb{N}\}$ is not finitely generated.
\end{ex}

\begin{lemma}
\emph{(\cite[Proposition 2.5]{BMS2013})\ }A semiring $S$ is $m$-Noetherian
if and only if every $S$-subsemimodule of a finitely generated free $S$%
-semimodule $S^{n}$ is finitely generated.
\end{lemma}

\begin{lemma}
\label{fg+N=fp}Let $S$ be an $m$-Noetherian semiring.

\begin{enumerate}
\item Every finitely generated $S$-semimodule is finitely presentable.

\item Every uniformly finitely generated $S$-semimodule is uniformly
finitely presented.
\end{enumerate}
\end{lemma}

\begin{Beweis}
\begin{enumerate}
\item This is \cite[Proposition 2.6]{BMS2013}.

\item This follows directly from the definitions.$\blacksquare $
\end{enumerate}
\end{Beweis}

\begin{lemma}
\label{inj-can}Let $\Lambda $ and $\Lambda ^{\prime }$ be two index sets. If 
$S$ is an $m$-Noetherian semiring, then the following canonical map is
injective%
\begin{equation*}
\beta _{(\Lambda ,\Lambda ^{\prime })}:S^{\Lambda }\otimes _{S}S^{\Lambda
^{\prime }}\longrightarrow S^{\Lambda \times \Lambda ^{\prime }},\text{ }%
f\otimes _{S}f^{\prime }\mapsto \lbrack (\lambda ,\lambda ^{\prime })\mapsto
f(\lambda )f^{\prime }(\lambda ^{\prime })].
\end{equation*}
\end{lemma}

\begin{Beweis}
Write $S^{\Lambda }=\underset{\longrightarrow }{\lim }M_{\lambda },$ a
direct limit of \emph{uniformly} finitely generated $S$-subsemimodules.
Since $S$ is $m$-Noetherian, $M_{\lambda }$ is uniformly finitely presented
for each $\lambda \in \Lambda .$ By Lemma \ref{fp-iso}, we have for each $%
\lambda \in \Lambda $ an isomorphism of $S$-semimodules $M_{\lambda }\otimes
_{S}S^{\Lambda ^{\prime }}\simeq M_{\lambda }^{\Lambda ^{\prime }}.$ It
follows that%
\begin{equation*}
S^{\Lambda }\otimes _{S}S^{\Lambda ^{\prime }}=\underset{\longrightarrow }{%
\lim }M_{\lambda }\otimes _{S}S^{\Lambda ^{\prime }}\simeq \underset{%
\longrightarrow }{\lim }(M_{\lambda }\otimes _{S}S^{\Lambda ^{\prime
}})\simeq \underset{\longrightarrow }{\lim }M_{\lambda }^{\Lambda ^{\prime
}}\hookrightarrow (S^{\Lambda })^{\Lambda ^{\prime }}\simeq S^{\Lambda
\times \Lambda ^{\prime }}.\blacksquare
\end{equation*}
\end{Beweis}

\begin{punto}
Let $(A,\mu ,\eta )$ be an $S$-semialgebra such that the canonical morphism
of $(A,A)$-bisemimodules%
\begin{equation}
\beta _{(A,A)}:R^{A}\otimes _{S}R^{A}\hookrightarrow R^{A\times A},\text{ }%
f\otimes _{S}f^{\prime }\mapsto \lbrack (a,a^{\prime })\mapsto f(a)f^{\prime
}(a^{\prime })]  \label{betaA}
\end{equation}%
is injective. Notice that the multiplication $\mu :A\otimes
_{S}A\longrightarrow A$ induces an (injective) $(A,A)$-bilinear map%
\begin{equation*}
\mu ^{\bullet }:S^{A}\longrightarrow S^{A\times A}.
\end{equation*}%
We define $A^{\circ }$ as the pullback in the following diagram in $_{A}%
\mathbb{S}_{A}:$%
\begin{equation*}
\xymatrix{A^{\circ} \ar@{^{(}.>}[d] \ar@{.>}[r] & S^A \otimes_S S^A
\ar@{^{(}->}[d]^{\beta_A} \\ A^* \ar[r]_(0.4){\mu ^{\bullet}} & S^{A \times
A}}
\end{equation*}%
equivalently%
\begin{eqnarray*}
A^{\circ } &=&\{(h,\sum\limits_{i=1}^{n}f_{i}\otimes _{S}g_{i})\in A^{\ast
}\oplus (S^{A}\otimes _{S}S^{A})\mid \mu ^{\bullet }(h)=\beta
_{(A,A)}(\sum\limits_{i=1}^{n}f_{i}\otimes _{S}g_{i})\} \\
&\simeq &\{h\in A^{\ast }\mid \exists \{(f_{i},g_{i})\}_{i=1}^{n}\subseteq
S^{A}\otimes _{S}S^{A}\text{ with }h(ab)=\sum%
\limits_{i=1}^{n}f_{i}(a)g_{i}(b)\text{ for all }a,b\in A\}.
\end{eqnarray*}
\end{punto}

\begin{lemma}
\label{f0}Let $A$ be an $S$-semialgebra and consider $A^{\ast }$ as an $%
(A,A) $-bisemimodule in the canonical way. If $S$ is $m$-Noetherian, then
the following are equivalent for $f\in A^{\ast }:$

\begin{enumerate}
\item $f\in A^{\circ };$

\item $Af$ is a finitely generated $S$-semimodule;

\item $\mu ^{\bullet }(f)\in \beta _{(A,A)}(A^{\circ }\otimes _{S}S^{A});$

\item $\mu ^{\bullet }(f)\in \beta _{(A,A)}(A^{\ast }\otimes _{S}S^{A});$

\item $fA$ is finitely generated;

\item $\mu ^{\bullet }(f)\in \beta _{(A,A)}(S^{A}\otimes _{S}A^{\circ });$

\item $\mu ^{\bullet }(f)\in \beta _{(A,A)}(S^{A}\otimes _{S}A^{\ast });$

\item $\mu ^{\bullet }(f)\in \beta _{(A,A)}(A^{\circ }\otimes _{S}S^{A})\cap
\beta _{(A,A)}(S^{A}\otimes _{S}A^{\circ }).$
\end{enumerate}
\end{lemma}

\begin{Beweis}
Notice that $A^{\circ }\hookrightarrow A^{\ast }$ is an $(A,A)$%
-subbisemimodule since it is -- by definition -- a pullback in $_{A}\mathbb{S%
}_{A}.$ The rest of the technical proof is now similar to that of \cite[%
Proposition 1.6]{AG-TW2000}.$\blacksquare $
\end{Beweis}

\begin{theorem}
\label{dual-co}Let $S$ be $m$-Noetherian and $A$ an $S$-semialgebra. If $%
A^{\circ }\hookrightarrow S^{A}$ is pure and $S^{A}$ is mono-flat, then $%
A^{\circ }$ is an $S$-semicoalgebra.
\end{theorem}

\begin{Beweis}
Since $S$ is $m$-Noetherian, the following $S$-linear map%
\begin{equation*}
\beta _{(A,A)}:S^{A}\otimes _{S}S^{A}\longrightarrow S^{A\times A},\text{ }%
f\otimes _{S}f^{\prime }\mapsto \lbrack (a,a^{\prime })\mapsto f(a)f^{\prime
}(a^{\prime })]
\end{equation*}%
is injective by Lemma \ref{inj-can}. Since $A^{\circ }\hookrightarrow S^{A}$
is pure and $S^{A}$ is mono-flat, we have the following canonical embeddings%
\begin{equation*}
A^{\circ }\otimes _{S}A^{\circ }\hookrightarrow S^{A}\otimes _{S}A^{\circ
}\hookrightarrow S^{A}\otimes _{S}S^{A}\overset{\beta _{(A,A)}}{%
\hookrightarrow }S^{A\times A}
\end{equation*}%
and the following canonical map%
\begin{equation*}
\widetilde{\beta }_{A}:A^{\circ }\otimes _{S}(A^{\circ }\otimes _{S}A^{\circ
})\hookrightarrow S^{A}\otimes _{S}(A^{\circ }\otimes _{S}A^{\circ
})\hookrightarrow S^{A}\otimes _{S}S^{A\times A}\overset{\beta _{(A,A\times
A)}}{\hookrightarrow }S^{A\times A\times A}.
\end{equation*}%
Moreover, for every $f\in A^{\circ },$ we have by Lemma \ref{f0}:%
\begin{equation*}
\mu ^{\ast }(f)\in \beta _{(A,A)}(A^{\circ }\otimes _{S}S^{A})\cap \beta
_{(A,A)}(S^{A}\otimes _{S}A^{\circ })=\beta _{(A,A)}(A^{\circ }\otimes
_{S}S^{A}\cap S^{A}\otimes _{S}A^{\circ })=\beta _{(A,A)}(A^{\circ }\otimes
_{S}A^{\circ }).
\end{equation*}%
Consider the following diagram%
\begin{equation*}
\xymatrix{S^A \ar[rrr]^{\mu ^{\bullet}} \ar[ddd]_{\mu ^{\bullet}} & & & S^{A
\times A} \ar[ddd]^{(id \times \mu)^{\bullet}}\\ & A^{\circ} \ar[ul]
\ar[r]^{\Delta} \ar[d]_{\Delta} & A^{\circ} \otimes_S A^{\circ} \ar[d]^{id
\otimes_S \Delta} \ar[ur] & \\ & A^{\circ} \otimes_S A^{\circ}
\ar[r]_(.4){\Delta \otimes_S id} \ar[dl] & A^{\circ} \otimes_S A^{\circ}
\otimes_S A^{\circ} \ar@{^{(}->}[dr]_{\tilde{\beta_(A,A)}} & \\ S^{A\times
A} \ar[rrr]_{(\mu \times id)^{\bullet}} & & & S^{A \times A \times A}}
\end{equation*}%
with%
\begin{equation*}
\mu ^{\bullet }:S^{A}\longrightarrow S^{A\times A},\text{ }f\mapsto \lbrack
(a,b)\mapsto f(ab)]\text{ and }\Delta =\mu _{\mid _{A^{\circ }}}^{\bullet
}:A^{\circ }\longrightarrow A^{\circ }\otimes _{S}A^{\circ }.
\end{equation*}%
Since the multiplication $\mu _{A}$ is associative, the outer rectangle is
clearly commutative. Moreover, all trapezoids are commutative. Since the
canonical amp $\widetilde{\beta }_{A}$ is injective, we conclude that the
inner diagram is commutative, \emph{i.e.} $\Delta $ is coassociative. It is
not difficult to show that the restriction of $\eta ^{\ast }:A^{\ast
}\longrightarrow S^{\ast }\simeq S$ is a counity for $A^{\circ
}.\blacksquare $
\end{Beweis}

\begin{theorem}
\label{dual-bi}Let $S$ be $m$-Noetherian.

\begin{enumerate}
\item If $B$ is an $S$-bisemialgebra with $B^{\circ }\hookrightarrow S^{B}$
pure and $S^{B}$ mono-flat, then $B^{\circ }$ is an $S$-bisemialgebra.

\item If $H$ is a Hopf $S$-semialgebra with $H^{\circ }\hookrightarrow S^{H}$
pure and $S^{H}$ mono-flat, then $H^{\circ }$ is a Hopf $S$-semialgebra.
\end{enumerate}
\end{theorem}

\end{document}